\documentclass[portuguese,a4paper,twoside,10pt]{article}
\usepackage[english]{babel}
\usepackage[T1]{fontenc}
\usepackage[latin1]{inputenc}
\usepackage{amssymb,amsmath,epsfig,amsfonts,euscript}
\usepackage{dsfont}
\usepackage[usenames,dvipsnames]{color}

\newtheorem{teo}{Theorem}[section]
\newtheorem{pro}[teo]{Proposition}
\newtheorem{cor}[teo]{Corollary}

\newtheorem{rem}[teo]{Remark}
\newtheorem{ex}{Example}

\newcounter{note}[section]
\newcounter{example}[section]

\newcommand{\er}{\mathbb{R}}

\newcommand{\dst}{\displaystyle}
\newcommand{\dem}{{\bf Proof }}
\newcommand{\fdem}{$\square$}

\newcommand{\sscr}{\scriptscriptstyle}

\newcommand{\nn}{\nonumber}
\newcommand{\mb}{\mathbf}

\newcommand{\jd}{$j=1,\hdots,d$}

\newcommand{\titulo}[1]{\begin{center}\mbox{} \\ \noindent \textit{\textbf{\Large #1}}\\\vspace{0.5cm}\end{center}}

\renewcommand{\abstract}[1]{{\small \noindent \textbf{Abstract:} #1\\}}

\newcommand{\keywords}[1]{{\small \noindent \textbf{Keywords:} #1\\}}

\begin{document}

\begin{center}
\titulo{Extremes of scale mixtures of multivariate time series}
\end{center}

\vspace{0.5cm}

\textbf{Helena Ferreira} Department of Mathematics, University of
Beira
Interior, Covilhã, Portugal \\(helena.ferreira@ubi.pt)\\

\textbf{Marta Ferreira} Center of Mathematics of Minho University, Braga, Portugal \\(msferreira@math.uminho.pt)\\

\abstract{Factor models have large potencial in the modeling of several natural and human phenomena. In this paper we consider a multivariate time series $\mb{Y}_n$, ${n\geq 1}$, rescaled through random factors $\mb{T}_n$, ${n\geq 1}$, extending some scale mixture models in the literature. We analyze its extremal behavior by deriving the maximum domain of attraction and the multivariate extremal index, which leads to new ways to construct multivariate extreme value distributions. The computation of the multivariate extremal index and the characterization of the tail dependence show the interesting property of these models that however much it is the dependence within and between factors  $\mb{T}_n$, ${n\geq 1}$, the extremal index of the model is unit whenever  $\mb{Y}_n$, ${n\geq 1}$, presents cross-sectional and sequencial tail independence. We illustrate with examples of thinned multivariate time series and multivariate autoregressive processes with random coefficients. An application of these latter to financial data is presented at the end.}


\section{Introduction}
Factor models have been used in the modeling of data within hydrology (Nadarajah \cite{nad06,nad09} 2006/2009, Nadarajah and Masoom \cite{nad+08} 2008), storm insurance (Lescourret and Robert, \cite{les+rob06} 2006), soil erosion in crops (Todorovic and Gani \cite{tod+gani87} 1987, Alpuim and Athayde \cite{alp+atha90} 1990), reliability (Alpuim and Athayde \cite{alp+atha90} 1990, Kotz \emph{et al.} \cite{kotz+2000} 2000), economy (Arnold, \cite{arn83} 1983) and finance (Ferreira and Canto e Castro, \cite{lccmf2} 2010).\\

Let $\mb{X}_n=(X_{n1}, \hdots, X_{nd})$, $n\geq 1$, be a $d$-variate sequence, such that $X_{nj}=Y_{nj}T_{nj}$, $j=1,\hdots, d$, where
\begin{itemize}
\item[(a)] $\mb{Y}=\{(Y_{n1}, \hdots, Y_{nd})\}_{n\geq 1}$ is a stationary sequence such that, $Y_{nj}$ has a Pareto-type distribution $F_{Y_j}$, $j=1,\hdots, d$, i.e., for each $j=1,\hdots, d$, there exists a positive constant  $\beta_j$ for which
    \begin{eqnarray}\label{Fy}
    F_{Y_j}(x)=1-x^{-\beta_j}l_{Y_j}(x),
    \end{eqnarray}
    with $l_{Y_j}$ a slowly varying function, i.e., $l_{Y_j}(ax)/l_{Y_j}(x)\to 1$, as $x\to\infty$, for all $a>0$,
\item[(b)] $\mb{T}=\{(T_{n1}, \hdots, T_{nd})\}_{n\geq 1}$ is a stationary sequence, independent of $\mb{Y}$, with support $\er_+^d$ and such that $E(T_{nj}^{\epsilon_j})<\infty$, for some $\epsilon_j>\beta_j$, $j=1,\hdots, d$.
\end{itemize}


This work is concerned with the extremal behavior of the multivariate time series $\mb{X}_n$, extending most of the factor models mentioned above. More precisely, we derive the max-domain of attraction (Section \ref{smda}), calculate the multivariate extremal index (Section \ref{smei}) and characterize the tail dependence (Section \ref{sti}).

The product $Y_{nj}T_{nj}$ can be seen as a random normalization of $Y_{nj}$ by $T_{nj}$, which is often required when modeling extremal behavior. For instance, if  $Y_{nj}$ is the rate of an extreme event and  $T_{nj}$ its average cost, then $Y_{nj}T_{nj}$ can be interpreted has the total cost of the extreme event. Products of two independent random variables where one of them is regularly varying have been addressed from both theoretical and applied points of view (Maulik \emph{et al.} \cite{mau+02} 2002, Lescourret and Robert \cite{les+rob06} 2006, Nadarajah \cite{nad06} 2006 and references therein).

Our motivation to the probabilistic study of extremes of multivariate sequences of products was originated from some particular models. Consider, for instance that $T_{nj}$ are Bernoulli distributed. Then $\mb{X}_n$ provides a model for multivariate data subjected to missing values. Extremes of univariate sequences with random missing values have been considered in Weissman and Cohen (\cite{weiss+co95}, 1995) as a particular case of some mixture models. Additional results on extremes of incomplete samples can be found in Mladenovic and Piterbarg (\cite{mla+pit06}, 2006) and Zhongquan and Wang (\cite{zhong+wang12}, 2012).

Li (\cite{li}, 2009) analyzed the tail dependence of the scale mixture $\mb{X}_n$ when $\mb{Y}_n=(Y_{n1},\hdots,Y_{nd})$ has multivariate extreme value distribution with standard Fr\'echet margins and $T_{nj}=T_n$, \jd.
Here we consider scale mixtures of multivariate sequences which are very flexible models for data exhibiting tail dependence and asymptotic tail independence such as, respectively, ARMAX and pARMAX sequences (Ferreira and Ferreira \cite{hf+mf6}).
We give particular emphasis to a model in which $\beta_j=\alpha/\gamma_j$, $\alpha,\gamma_j>0$, $j=1,\hdots, d$, generalizing the results of Lescourret and Robert (\cite{les+rob06}, 2006) (Section \ref{spc}). An application to financial data will be provided at the end (Section \ref{saplic}).

\section{Preliminary results and max-domain of attraction}\label{smda}

We start with some properties of $\{\mb{X}_n\}_n\geq 1$, that will be used along the paper. We use notation $r_j=E(T_{nj}^{\beta_j})$ along the paper.

\begin{pro}\label{pFx}
For each $j=1,\hdots, d$, $\{X_{nj}\}_{n\geq 1}$ is a stationary sequence having Pareto-type distribution.
\end{pro}
\dem First, observe that
$$
\dst\lim_{x\to\infty}\frac{P(Y_{nj}T_{nj}>x)}{P(Y_{nj}>x)}=
\lim_{x\to\infty}\frac{\dst \int x^{-\beta_j}z^{\beta_j}l_{Y_j}(x/z)dP_{T_{nj}}(z)}{x^{-\beta_j}l_{Y_j}(x)}=r_j,
$$
where the last step is due to the dominated convergence theorem and by using the Potter bounds of regularly varying functions (Bingham \emph{et al.}, \cite{bingham+} 1987; Theorem 1.5.6.). Therefore, for large $x$,
\begin{eqnarray}\label{Fx}
1-F_{X_j}(x)=P(X_{nj}>x)=x^{-\beta_j}l_{Y_j}(x)r_j(1+o(1)):=x^{-\beta_j}l_{X_j}(x),
 \end{eqnarray}
where it is immediately seen that $l_{X_j}$ is a slowly varying function. \fdem\\

In the sequel we  denote $U_{X_j}{(x)}$ and $U_{Y_j}{(x)}$ the quantile functions, $F_{X_j}^{-1}(1-1/x)$ and $F_{Y_j}^{-1}(1-1/x)$, respectively.

Given (\ref{Fy}), we can state
\begin{eqnarray}\label{F-1y}
U_{Y_j}(x)=x^{1/\beta_j}l_{U_{Y_j}}(x),
\end{eqnarray}
where  $l_{U_{Y_j}}$ is a slowly varying function, and by Proposition \ref{pFx}, we can also write
\begin{eqnarray}\nn\label{F-1x}
U_{X_j}(x)=x^{1/\beta_j}l_{U_{X_j}}(x),
\end{eqnarray}
where  $l_{U_{X_j}}$ is a  slowly varying function. By the Bruyn conjugate concept (Beirlant \emph{et al.} \cite{beirl+} 2004, Proposition 2.5), we have that, for large $x$,
\begin{eqnarray}\nn
U_{X_j}(x)=x^{1/\beta_j}l_{X_j}^{1/\beta_j}(x^{1/\beta_j}l_{U_{X_j}}(x))(1+o(1)).
\end{eqnarray}
By (\ref{Fx}) and (\ref{F-1y}), we have  for large $x$
\begin{eqnarray}\label{F-1xy}
U_{X_j}(x)=x^{1/\beta_j}l_{U_{Y_j}}(x)r_j^{1/\beta_j}(1+o(1))=U_{Y_j}\left(r_j\,x\right)(1+o(1)).
\end{eqnarray}

\begin{pro}\label{putc}
The upper tail copula function of $\mb{X}$ is given by
\begin{eqnarray}\nn
\Lambda_{\mb{X}}(x_1,...,x_d)=E\left(\Lambda_{\mb{Y}}\left(\frac{T_1^{\beta_1}x_1}{r_1},\hdots,\frac{T_d^{\beta_d}x_d}{r_d}\right)\right),
\end{eqnarray}
with $(x_1,...,x_d)\in\overline{\er}_+^d=[0,\infty]^d\backslash\{(\infty,\hdots,\infty)\}$, where $\mb{T}=(T_1,\hdots,T_d)$ is a random vector distributed as $\mb{T}_n=(T_{n1},\hdots,T_{nd})$ and provided that the upper tail copula function of $\mb{Y}_n$ exists, i.e., the limit
\begin{eqnarray}\label{utcopx}
\dst\Lambda_{\mb{Y}}(x_1,...,x_d)=\lim_{t\to\infty}tP\left(\bigcap_{j=1}^d\left\{Y_{1j}>U_{Y_j}\left(t/x_j\right)\right\}\right)
\end{eqnarray}
is finite.
\end{pro}
\dem The upper tail copula function of $\mb{X}$ is defined by
\begin{eqnarray}\nn
\dst\Lambda_{\mb{X}}(x_1,...,x_d)=\lim_{t\to\infty}tP\left(\bigcap_{j=1}^d\left\{X_{1j}>U_{X_j}\left(t/x_j\right)\right\}\right)
\end{eqnarray}

The result follows immediately by applying (\ref{F-1xy}) and the dominated convergence theorem, since
\begin{eqnarray}\label{utcopxy}
\begin{array}{rl}
&\dst \lim_{t\to\infty}tP\left(\bigcap_{j=1}^d\left\{X_{1j}>U_{X_j}\left(t/x_j\right)\right\}\right)\vspace{0.35cm}\\
=&\dst\lim_{t\to\infty} t\int P\left(\bigcap_{j=1}^d\left\{Y_{1j}>U_{Y_j}\left(\frac{r_j\,t}{z_j^{\beta_j}\,x_j}\right)\right\}\right)dP_{(T_{11},\hdots,T_{1d})}(z_1,\hdots,z_d)\vspace{0.35cm}\\
=&\dst \int \Lambda_{\mb{Y}}\left(\frac{z_1^{\beta_1}x_1}{r_1},\hdots,\frac{z_d^{\beta_d}x_d}{r_d}\right)dP_{(T_{11},\hdots,T_{1d})}(z_1,\hdots,z_d)
\vspace{0.35cm}\\
=&\dst E\left(\Lambda_{\mb{Y}}\left(\frac{T_1^{\beta_1}x_1}{r_1},\hdots,\frac{T_d^{\beta_d}x_d}{r_d}\right)\right).\,\square\\
\end{array}
\end{eqnarray}
\smallskip

\begin{rem}
If $\beta_j=\beta$ and $T_{nj}=T$, \jd, then $\Lambda_{\mb{X}}=\Lambda_{\mb{Y}}$, corresponding to the case considered in Li (\cite{li}, 2009). If $\Lambda_{\mb{Y}}(x_1,...,x_d)=\sum_{i=1}^dx_j$, we also obtain $\Lambda_{\mb{X}}=\Lambda_{\mb{Y}}$, for any choice of $\mb{T}$. In the last section we give particular attention to the case $\Lambda_{\mb{Y}}(x_1,...,x_d)=\bigwedge_{i=1}^dx_j$.
\end{rem}
Normalized levels $u_{nj}^{(\tau_j)}$ of $X_j$ are such that $n(1-F_{X_j}(u_{nj}^{(\tau_j)}))\to\tau_j>0$, as $n\to\infty$, i.e., $u_{nj}^{(\tau_j)}=U_{X_j}\left(n/\tau_j\right)(1+o(1))$, as $n\to\infty$. Thus, by (\ref{F-1xy}) and for large $n$, $X_j$ has normalized levels
$$
u_{nj}^{(\tau_j)}=U_{Y_j}\left(nr_j/\tau_j\right)(1+o(1)).
$$
Consider $\mb{u}_n(\boldsymbol{x})=(u_{n1}^{(\tau_1(x_1))},\hdots,u_{nd}^{(\tau_d(x_d))})
=(nx_1,\hdots,nx_d)$ a vector of normalized levels of $\mb{X}$ and, for each \jd, let $\{\widehat{X}_{nj}\}_{n\geq 1}$ be an i.i.d.~sequence with the same marginal distribution as $\{X_{nj}\}_{n\geq 1}$.
Denote $\widehat{\mb{M}}_{n}=(\widehat{{M}}_{n1}, \hdots,\widehat{{M}}_{nd})$ the vector of  the componentwise maxima $\widehat{{M}}_{nj}=\bigvee_{i=1}^n\widehat{X}_{ij}$, $j=1,\hdots,d$.
\begin{pro}\label{pmev}
We have $F_{\mb{X}_1}$  in the domain of attraction of $G_{\mb{X}}$, that is,
$$
\dst \lim_{n\to\infty}P\left(\widehat{\mb{M}}_{n}\leq \mb{u}_n(\boldsymbol{x})\right)=G_{\mb{X}}(x_1,...,x_d),
$$
with
\begin{eqnarray}\nn
G_{\mb{X}}(x_1,...,x_d)=\exp\left\{-E\left(-\log G_{\mb{Y}}\left(\frac{x_1r_1}{T_1^{\beta_1}},\hdots,\frac{x_dr_d}{T_d^{\beta_d}}\right)\right)\right\},
\end{eqnarray}
provided that $F_{\mb{Y}_1}$ is in the domain of attraction of $G_{\mb{Y}}(x_1,...,x_d)$, where both $G_{\mb{X}}$ and $G_{\mb{Y}}$ have unit Fr\'echet marginals.
\end{pro}
\dem
Just observe that
$$
\dst \lim_{n\to\infty}P\left(\bigcap_{j=1}^d\left\{\widehat{M}_{nj}\leq U_{X_j}\left(nx_j\right)\right\}\right)=\dst \lim_{n\to\infty}\exp\left\{-nP\left(\bigcup_{j=1}^d\left\{X_{1j}>U_{X_j}\left(nx_j\right)\right\}\right)\right\}.\,\square
\nn\\\nn $$
Now the proof runs along the same lines as in  (\ref{utcopxy}). \fdem\\

In the following, for any vector $\mb{z}$ and $A\subset\{1,\hdots,d\}$, $\mb{z}_A$ denotes the sub-vector of $\mb{z}$ with indices in $A$.

\begin{ex}
If $(T_1,\hdots,T_d)$ has support in $\{0,1\}^d$ and $P\left(\bigcap_{i\in J}T_i=1,\bigcap_{i\in D\backslash J}T_i=0\right)=p(J)$, $\emptyset\not =J\subset D=\{1,\hdots,d\}$, then $G_{\mb{X}}$ corresponds to a geometric mean of the marginal distributions of $G_{\mb{Y}}$. In this case, we have
\begin{eqnarray}\nn
\begin{array}{rl}
\dst G_{\mb{X}}(x_1,\hdots,x_d)=&\dst \prod_{\emptyset\not=J\subset D}G_{\mb{Y}_J}^{p(J)}(r_1x_1,\hdots,r_dx_d)_{J}
= \dst \prod_{\emptyset\not=J\subset D}G_{\mb{Y}_J}\left(\frac{r_1x_1}{p(J)},\hdots,\frac{r_dx_d}{p(J)}\right)_{J}
\vspace{0.35cm}\\
=& \dst \prod_{\emptyset\not=J\subset D}G_{\mb{Y}_J}\left(\frac{p(\{1\})\,x_1}{p(J)},\hdots,\frac{p(\{d\})\,x_d}{p(J)}\right)_{J}.
\end{array}
\end{eqnarray}
We illustrate the result with some choices for $G_{\mb{Y}}$.

If the stationary sequence $\mb{Y}$ has common copula logistic, i.e.,
$$
C_{\mb{Y}_n}(u_1,\hdots,u_d)=\exp\left\{-\left(\sum_{j=1}^d(-\log u_j)^{-1/\alpha}\right)^{\alpha}\right\}
$$
then we obtain for $G_{\mb{Y}}$ the logistic distribution and
$$
G_{\mb{X}}(x_1,\hdots,x_d)=\exp\left\{-\sum_{\emptyset\not=J\subset D}\left(\sum_{i\in J}(\beta_{Ji}x_i)^{-1/\alpha}\right)^{\alpha}\right\},
$$
with $\beta_{Ji}={p(\{i\})}/{p(J)}$, which is an asymmetric logistic distribution already found in Tawn (\cite{tawn90}, 1990), by following a different probabilistic approach. The parameters $\beta_{Ji}$ increase the variability within the tail dependence coefficients regarding the departure distribution $G_{\mb{Y}}$. For $\mb{X}=(X_1,\hdots,X_d)$ with distribution $G_{\mb{X}}$, we have, for instance,
\begin{eqnarray}\nn
\begin{array}{rl}
\dst \Lambda_{(X_i,X_j)}(1,1)=&\dst \lim_{t\to\infty}
P(F_{X_i}(X_i)>1-1/t|F_{X_j}(X_j)>1-1/t)
=2-\sum_{J\in\mathcal{F}_{\{i,j\}}}
\left(\beta_{Ji}^{-1/\alpha}+\beta_{Jj}^{-1/\alpha}\right)^{\alpha},
\end{array}
\end{eqnarray}
where $\mathcal{F}_{\{i,j\}}$ is the family of subsets of $D$ containing $\{i,j\}$, expression that presents a larger number of possibilities than the value $\Lambda_{(Y_i,Y_j)}(1,1)=2-2^{\alpha}$ of the symmetric logistic distribution. Note that $\sum_{J\in \mathcal{F}_{\{i\}}}\beta_{Ji}=1$.

Now, suppose that  $G_{\mb{Y}}(x_1,\hdots,x_d)=\prod_{l=1}^{\infty}
\prod_{k=-\infty}^\infty\bigwedge_{j=1}^d\exp\left(-a_{lkj}x_j^{-1}\right)$,
where $\{a_{lkj},l\geq 1, -\infty<k<\infty,1\leq j\leq d\}$ are real non negative constants satisfying $\sum_{l=1}^{\infty}
\sum_{k=-\infty}^\infty a_{lkj}=1$, \jd. This is the attractor MEV when $\mb{Y}$ is an M4 process (Smith and Weissman \cite{smith+weissman}, 1996) and we have
$$
G_{\mb{X}}(x_1,\hdots,x_d)=\prod_{l=1}^{\infty}
\prod_{k=-\infty}^\infty\prod_{\emptyset\not =J\subset D}\bigwedge_{j\in J}\exp\left(-\frac{a_{lkj}}{\beta_{Jj}}x_j^{-1}\right)
$$
and the bivariate tail dependence summarized by
\begin{eqnarray}\nn
\begin{array}{rl}
\dst \Lambda_{(X_i,X_j)}(1,1)
=2-\sum_{l=1}^{\infty}
\sum_{k=-\infty}^\infty\sum_{J\in\mathcal{F}_{\{i,j\}}}
\left(\frac{a_{lki}}{\beta_{Ji}}\vee \frac{a_{lkj}}{\beta_{Jj}}\right)\leq 2-\sum_{l=1}^{\infty}
\sum_{k=-\infty}^\infty
\left({a_{lki}}\vee {a_{lkj}}\right)=\Lambda_{(Y_i,Y_j)}(1,1).
\end{array}
\end{eqnarray}
\end{ex}

\section{The multivariate extremal index}\label{smei}

The extremal index measures the tendency of clusters occurrence, a phenomena commonly observed in real data. In this section we will compute the multivariate extremal index of $\mb{X}_n$ (Nandagopalan \cite{nand}, 1990). We start by analyzing some long range and local dependence conditions that will make easier its calculation.

\begin{pro}
If $\mb{Y}$ and $\mb{T}$ are strong-mixing, then $\mb{X}$ is strong-mixing.
\end{pro}
\dem Consider the events $A$ and $B$, respectively, in $\sigma$-algebras generated by $\{\mb{X}_1,\hdots,\mb{X}_p\}$ and $\{\mb{X}_{p+s},\mb{X}_{p+s+1},\hdots\}$, i.e., $A\in\sigma(\mb{X}_1,\hdots,\mb{X}_p)$ and $B\in\sigma(\mb{X}_{p+s},\hdots)$.
Given the independence between $\mb{Y}$ and $\mb{T}$, we can state
$$
\begin{array}{rl}
&P(A\cap B)=E(P(A'\cap B'|\mb{Y}))= E\left(P(A'|\mb{Y}_1,\hdots,\mb{Y}_p)P(B'|\mb{Y}_{p+s},\hdots)\right)+O(\alpha_{\mb{T}}(s)),
\end{array}\nn\\\nn
$$
where $A'\in\sigma(\mb{Y}_1,\hdots,\mb{Y}_p)$,  $B'\in\sigma(\mb{Y}_{p+s},\hdots)$ and $\alpha_{\mb{T}}(s)$ denotes the mixture coefficients of the sequence  $\mb{T}$. The result follows by Doukhan \cite{douk95} (1995, Theorem 3 in Section 1.2.2), since
$$
E\left(P(A'|\mb{Y}_1,\hdots,\mb{Y}_p)P(B'|\mb{Y}_{p+s},\hdots)\right)
\leq E\left(P(A'|\mb{Y}_1,\hdots,\mb{Y}_p)\right)E\left(P(B'|\mb{Y}_{p+s},\hdots)\right)
+4\alpha_{\mb{Y}}(s),\\
$$
where $\alpha_{\mb{Y}}(s)$ denotes the mixture coefficients of the sequence  $\mb{Y}$. \fdem\\

Consider $\mb{u}_n(\boldsymbol{\tau})=(u_{n1}^{(\tau_1)},\hdots,u_{nd}^{(\tau_d)})$ a vector of normalized levels of $\mb{X}$ and ${\mb{M}_{1,n}\equiv} \mb{M}_{n}=({M}_{n1}, \hdots,{M}_{nd})$ the vector of  the componentwise maxima ${M}_{nj}=\bigvee_{i=1}^nX_{ij}$, $j=1,\hdots,d$. As a consequence of the strong-mixing condition, we have
$$
\begin{array}{rl}
\dst\lim_{n\to\infty}P\left(\mb{M}_n\leq \mb{u}_n^{(\boldsymbol{\tau})}\right) =
&\dst\lim_{n\to\infty}P\left({M}_{n1}\leq {u}_{n1}^{({\tau_1})}, \hdots,{M}_{nd}\leq {u}_{nd}^{({\tau_d})} \right)\vspace{0.35cm}\\
=&\exp\left\{-\dst\lim_{n\to\infty}nP\left(\mb{X}_1\not\leq \mb{u}_n^{(\boldsymbol{\tau})},\bigcap_{i=2}^{[n/k_n]}\mb{X}_{i}\leq \mb{u}_n^{(\boldsymbol{\tau})}\right)\right\}
\end{array}\nn\\\nn
$$
where $k_n\to\infty$ and $n/k_n\to\infty$, as $n\to\infty$. Therefore, $\mb{X}$ has multivariate extremal index
$$
\theta_{\mb{X}}(\tau_1,\hdots,\tau_d)=\frac{-\log \lim_{n\to\infty}P\left(\mb{M}_n\leq \mb{u}_n^{(\boldsymbol{\tau})}\right)}{-\log \lim_{n\to\infty}P\left(\widehat{\mb{M}}_n\leq \mb{u}_n^{(\boldsymbol{\tau})}\right)},
$$
if and only if $\lim_{n\to\infty}nP\left(\mb{X}_1\not\leq \mb{u}_n^{(\boldsymbol{\tau})},{ \mb{M}_{2,[n/k_n]}}\leq \mb{u}_n^{(\boldsymbol{\tau})}\right)$ exists and, if so, we have
\begin{eqnarray}\label{mindext}
\theta_{\mb{X}}(\tau_1,\hdots,\tau_d)=\dst\lim_{n\to\infty}\frac{P\left(\mb{X}_1\not\leq \mb{u}_n^{(\boldsymbol{\tau})},\mb{M}_{ 2,[n/k_n]}\leq \mb{u}_n^{(\boldsymbol{\tau})}\right)}{P\left(\mb{X}_1\not\leq \mb{u}_n^{(\boldsymbol{\tau})}\right)}.
\end{eqnarray}
This function relates the two multivariate extreme value distributions arising from the maxima in $\{\mb{X}_n\}_{n\geq 1}$ and $\{\widehat{\mb{X}}_n\}_{n\geq 1}$. We analyze its values under additional assumptions on the scale sequence $\mb{Y}$.\\

Let
$$
\frac{\mb{r}\,n}{\boldsymbol{\tau}\,\mb{T}_i^{\boldsymbol{\beta}}}
=\left(\frac{r_1\,n}{\tau_1\,T_{i1}^{\beta_1}},\hdots, \frac{r_d\,n}{\tau_d\,T_{id}^{\beta_d}}\right)
$$
and $\mb{U}_{\mb{Y}_i}(\mb{x})=(U_{Y_{i1}}(x_1),\hdots, U_{Y_{id}}(x_d))$. About the numerator in (\ref{mindext}) we can write
\begin{eqnarray}\label{mindext1}
\begin{array}{rl}
&\dst nP\left(\mb{X}_1\not\leq \mb{u}_n^{(\boldsymbol{\tau})},{ \mb{M}_{2,[n/k_n]}}\leq \mb{u}_n^{(\boldsymbol{\tau})}\right)\vspace{0.35cm}\\
=&\dst
nP\left(\mb{Y}_1\not\leq { \mb{U}_{\mb{Y}_1}}\left(\frac{\mb{r}\,n}{\boldsymbol{\tau}\,\mb{T}_1^{\boldsymbol{\beta}}}\right),\bigcap_{i=2}^{[n/k_n]}\mb{Y}_i\leq { \mb{U}_{\mb{Y}_i}}\left(\frac{\mb{r}\,n}{\boldsymbol{\tau}\,\mb{T}_i^{\boldsymbol{\beta}}}\right)\right).
\vspace{0.35cm}\\
\end{array}
\end{eqnarray}
In some cases, the sequence $\mb{Y}$ is such that
$$
\begin{array}{rl}
&\dst \lim_{n\to\infty} nP\left(\mb{Y}_1\not\leq \mb{U}_{\mb{Y}_i}\left(\frac{n}{\boldsymbol{\tau}}\right),\bigcap_{i=2}^{[n/k_n]}\mb{Y}_i\leq \mb{U}_{\mb{Y}_i}\left(\frac{n}{\boldsymbol{\tau}}\right)\right)\vspace{0.35cm}\\
=&\dst
\dst \lim_{n\to\infty} nP\left(\mb{Y}_1\not\leq \mb{U}_{\mb{Y}_i}\left(\frac{n}{\boldsymbol{\tau}}\right),\bigcap_{i=2}^{k}\mb{Y}_i\leq \mb{U}_{\mb{Y}_i}\left(\frac{n}{\boldsymbol{\tau}}\right)\right),
\end{array}
$$
for some levels $\mb{u}_n^{(\boldsymbol{\tau})}=\mb{U}_{\mb{Y}}\left(\frac{n}{\boldsymbol{\tau}}\right)$ and for some finite $k$, which makes easier the calculation of the multivariate extremal index. Situations leading to this easier procedure are highlighted in the next results.

\begin{pro}\label{pdk}
Consider $\mb{a}=(a_1,\hdots,a_d)$, $\mb{b}=(b_1,\hdots,b_d)$, $T_{nj}$ with support in $[a_j,b_j]$, $j=1,\hdots,d$, and $\mb{v}_n^{(\boldsymbol{\tau}^*(\mb{z}))}=\left(U_{Y_1}
\left(\frac{r_1n}{z_1^{\beta_1}\tau_1}\right),\hdots,U_{Y_d}
\left(\frac{r_dn}{z_d^{\beta_d}\tau_d}\right)\right)$, for all $z_j\in[a_j,b_j]$, $j=1,\hdots,d$. If $\mb{Y}$ satisfies condition ${\tt D}^{(k)}\left(\mb{v}_n^{(\boldsymbol{\tau}^*(\mb{a}))},\mb{v}_n^{(\boldsymbol{\tau}^*(\mb{b}))}\right)$, defined by
\begin{eqnarray}\label{pdk1}
\lim_{n\to\infty} n\sum_{i=k}^{[n/k_n]}P\left(\mb{Y}_{1}\not\leq \mb{v}_n^{(\boldsymbol{\tau}^*(\mb{b}))},\mb{Y}_{i}\leq \mb{v}_n^{(\boldsymbol{\tau}^*(\mb{a}))},\mb{Y}_{i+1}\not\leq \mb{v}_n^{(\boldsymbol{\tau}^*(\mb{b}))}\right)=0,
\end{eqnarray}
then ${\tt D}^{(k)}\left(\mb{u}_n^{(\boldsymbol{\tau})}\right)$, defined by
\begin{eqnarray}\nn
\lim_{n\to\infty} n\sum_{i=k}^{[n/k_n]}P\left(\mb{X}_{1}\not\leq \mb{u}_n^{(\boldsymbol{\tau})},\mb{X}_{i}\leq \mb{u}_n^{(\boldsymbol{\tau})},\mb{X}_{i+1}\not\leq \mb{u}_n^{(\boldsymbol{\tau})}\right)=0,
\end{eqnarray}
holds for $\mb{X}$, $\boldsymbol{\tau}=(\tau_1,\hdots,\tau_d)$.
\end{pro}
\dem Observe that, if we take  $\mb{u}_n^{(\boldsymbol{\tau})}=\left(U_{X_1}\left(n/\tau_1 \right),\hdots,U_{X_d}\left(n/\tau_d\right)\right)$, we have
\begin{eqnarray}\nn
\begin{array}{rl}
&\dst \lim_{n\to\infty} n\sum_{i=k}^{[n/k_n]}P\left(\mb{X}_{1}\not\leq \mb{u}_n^{(\boldsymbol{\tau})},\mb{X}_{i}\leq \mb{u}_n^{(\boldsymbol{\tau})},\mb{X}_{i+1}\not\leq \mb{u}_n^{(\boldsymbol{\tau})}\right)\vspace{0.35cm}\\
=&\dst \lim_{n\to\infty}\int n\sum_{i=k}^{[n/k_n]}P\left(\mb{Y}_{1}\not\leq \mb{v}_n^{(\boldsymbol{\tau}^*(\mb{z}_1))},\mb{Y}_{i}\leq \mb{v}_n^{(\boldsymbol{\tau}^*(\mb{z}_{i}))},\mb{Y}_{i+1}\not\leq \mb{v}_n^{(\boldsymbol{\tau}^*(\mb{z}_{i+1}))}\right)
dP_{(\mb{T}_{1},\mb{T}_{i},\mb{T}_{i+1})}(\mb{z}_1,\mb{z}_{i},\mb{z}_{i+1})
\end{array}
\end{eqnarray}

Now, just observe that
$$
\begin{array}{rl}
&P\left(\mb{Y}_{1}\not\leq \mb{v}_n^{(\boldsymbol{\tau}^*(\mb{z}_1))},\mb{Y}_{i}\leq \mb{v}_n^{(\boldsymbol{\tau}^*(\mb{z}_{i}))},\mb{Y}_{i+1}\not\leq \mb{v}_n^{(\boldsymbol{\tau}^*(\mb{z}_{i+1}))}\right)\vspace{0.35cm}\\
\leq &P\left(\mb{Y}_{1}\not\leq \mb{v}_n^{(\boldsymbol{\tau}^*(\mb{b}))},\mb{Y}_{i}\leq \mb{v}_n^{(\boldsymbol{\tau}^*(\mb{a}))},\mb{Y}_{i+1}\not\leq \mb{v}_n^{(\boldsymbol{\tau}^*(\mb{b}))}\right).\,\square
\end{array}\nn\\\nn
$$
\smallskip

Under the condition (\ref{pdk1}), smooth oscillations around $\mb{v}_n^{(\boldsymbol{\tau}^*(\mb{b}))}$ by $\mb{Y}$ (in the sense that its values which are at least lag-$k$ apart no longer exceed $\mb{v}_n^{(\boldsymbol{\tau}^*(\mb{a}))}>\mb{v}_n^{(\boldsymbol{\tau}^*(\mb{b}))}$) are not followed by subsequent exceedances of $\mb{v}_n^{(\boldsymbol{\tau}^*(\mb{b}))}$. In the Example \ref{emrarmax}, we provide an illustration of such local behavior which leads to the condition ${\tt D}^{(k)}\left(\mb{u}_n^{(\boldsymbol{\tau})}\right)$ for the scale model $\mb{X}$.

In order to avoid restrictions on the support of $T_{nj}$, \jd, we propose now a greater restriction to the oscillations of $\mb{Y}$ around normalized levels. First observe that
$$
P\left(\mb{Y}_{1}\not\leq \mb{v}_n^{(\boldsymbol{\tau}^*(\mb{z}_1))},\mb{Y}_{i}\leq \mb{v}_n^{(\boldsymbol{\tau}^*(\mb{z}_{i}))},\mb{Y}_{i+1}\not\leq \mb{v}_n^{(\boldsymbol{\tau}^*(\mb{z}_{i+1}))}\right)
\leq P\left(\mb{Y}_{1}\not\leq \mb{v}_n^{(\boldsymbol{\tau}^*(\mb{z}))},\mb{Y}_{i+1}\not\leq \mb{v}_n^{(\boldsymbol{\tau}^*(\mb{z}))}\right)
$$
with
$$
\mb{v}_n^{(\boldsymbol{\tau}^*(\mb{z}))}=
\mb{U}_{\mb{Y}_1}\left(\frac{r_1 n}{(z_{11}\vee z_{i+1,1})^{\beta_1}\tau_1},\hdots,\frac{r_d n}{(z_{1d}\vee z_{i+1,d})^{\beta_d}\tau_d}\right).
$$
Therefore, we can take for $\mb{Y}$ a sequence satisfying the following restriction on the local occurrence of exceedances of $\mb{v}_n^{(\boldsymbol{\tau})}$ which are lag-$k$ apart:
\begin{eqnarray}\label{kdep}
\dst\lim_{n\to\infty}n\sum_{i=k}^{[n/k]} P\left(\mb{Y}_{1}\not\leq \mb{v}_n^{(\boldsymbol{\tau})},\mb{Y}_{i+1}\not\leq \mb{v}_n^{(\boldsymbol{\tau})}\right)=0.
\end{eqnarray}
This condition limits the size of the clusters of high levels and is satisfied by $k$-dependent sequences  $\mb{Y}$, i.e., $\mb{Y}_{n}$ and $\mb{Y}_{m}$ are independent whenever $|n-m|>k$.

\begin{pro}
If $\mb{Y}$ satisfies (\ref{kdep}) for all $\mb{v}_n^{(\boldsymbol{\tau})}=\mb{U}_{\mb{Y}_1}\left(\frac{n}{\boldsymbol{\tau}}\right)$, then $\mb{X}$ satisfies ${\tt D}^{(k)}\left(\mb{u}_n^{(\boldsymbol{\tau})}\right)$, for any sequence of normalized levels $\mb{u}_n^{(\boldsymbol{\tau})}=\mb{U}_{\mb{X}_1}\left(\frac{n}{\boldsymbol{\tau}}\right)$.
\end{pro}

\begin{cor}
If $\mb{Y}$ is $k$-dependent, then $\mb{X}$ satisfies  ${\tt D}^{(k)}\left(\mb{u}_n^{(\boldsymbol{\tau})}\right)$, for any sequence of normalized levels $\mb{u}_n^{(\boldsymbol{\tau})}$.
\end{cor}

If $\mb{X}$ is strong-mixing and satisfies ${\tt D}^{(k)}\left(\mb{u}_n^{(\boldsymbol{\tau})}\right)$ condition, for all $\boldsymbol{\tau}\in \er_+^d$, then it has multivariate extremal index given by
\begin{eqnarray}\label{mindextdk}
\theta_{\mb{X}}(\tau_1,\hdots,\tau_d)=\lim_{n\to\infty}\frac{P\left(\mb{X}_1\leq \mb{u}_n^{(\boldsymbol{\tau})},\hdots,\mb{X}_{k-1}\leq \mb{u}_n^{(\boldsymbol{\tau})},\mb{X}_{k}\not\leq \mb{u}_n^{(\boldsymbol{\tau})}\right)}{P\left(\mb{X}_1\not\leq \mb{u}_n^{(\boldsymbol{\tau})}\right)},
\end{eqnarray}
whenever this limit exists (Chernick \emph{et al.} \cite{chern+91} 1991, Ferreira \cite{hf94} 1994) . In this case,
\begin{eqnarray}\label{mevdk}
\begin{array}{c}
\dst\lim_{n\to\infty}P\left(\mb{M}_n\leq \mb{u}_n^{(\boldsymbol{\tau})}\right) =\exp\left\{-\dst\lim_{n\to\infty}nP\left(\mb{X}_1\leq \mb{u}_n^{(\boldsymbol{\tau})},\hdots, \mb{X}_{k-1}\leq \mb{u}_n^{(\boldsymbol{\tau})} ,\mb{X}_{k}\not\leq \mb{u}_n^{(\boldsymbol{\tau})}\right)\right\}.
\end{array}
\end{eqnarray}

The probability in the second term can be computed throughout the tail behavior of $(\mb{Y}_1,\hdots,\mb{Y}_k)$. In order to improve such approach, define
\begin{eqnarray}\label{LambdaY1k}
\Lambda_{(\mb{Y}_1,\hdots,\mb{Y}_k)}(\mb{x}_1,\hdots,\mb{x}_k)=\dst\lim_{t\to\infty}
t\,P\left({ \mb{Y}_1\not\leq \mb{U}_{\mb{Y}_1}\left(\frac{t}{\mb{x}_1}\right),\hdots,\mb{Y}_k\not\leq \mb{U}_{\mb{Y}_k}\left(\frac{t}{\mb{x}_k}\right)} \right),
\end{eqnarray}
with ${ \mb{x}_j\in\overline{\er}_+^{\,d}}$, $j=1,\hdots,k$. This definition extends the upper tail dependence concept of Definition 2.3 in Ferreira and Ferreira (\cite{hf+mf2}, 2012).\\

Note that, for $\mb{x}_i$ with null components $x_{ij}$ with $j\in\{1,\hdots,d\}\backslash J_i$, $J_i\not=\emptyset$, $i=1,\hdots,k$, we have
$$
\Lambda_{(\mb{Y}_1,\hdots,\mb{Y}_k)}(\mb{x}_1,\hdots,\mb{x}_k)
=\Lambda_{((\mb{Y}_1)_{J_1},\hdots,(\mb{Y}_k)_{J_k})}
((\mb{x}_1)_{J_1},\hdots,(\mb{x}_k)_{J_k}),
$$
where $(\mb{x}_i)_{J_i}$ denotes the sub-vector of $\mb{x}_i$ having components indexed in $J_i$. On the other hand, if $\mb{x}$ has some null component, we have $\Lambda_{\mb{Y}_1}(\mb{x})=0$.\\

We are going to apply the result in (\ref{mevdk}) to our model and derive an expression for the multivariate extremal index, which will depend on the upper tail dependence function of $(\mb{Y}_1,\hdots,\mb{Y}_k)$ given in (\ref{LambdaY1k}) and also on the dependence between random vectors $\mb{T}_1,\hdots,\mb{T}_k$.

\begin{pro}
If (\ref{LambdaY1k}) holds for all $\mb{x}_1,\hdots,\mb{x}_k\in { \overline{\er}_+^{\,d}}$ then
\begin{eqnarray}\nn
\begin{array}{rl}
&\dst\lim_{n\to\infty} n\, P\left(\mb{X}_1\leq \mb{u}_n^{(\boldsymbol{\tau})},\hdots,\mb{X}_{k-1}\leq \mb{u}_n^{(\boldsymbol{\tau})},\mb{X}_{k}\not\leq \mb{u}_n^{(\boldsymbol{\tau})}\right)\vspace{0.35cm}\\
=& \dst E\left(\sum_{\emptyset\subseteq I\subset\{1,\hdots,k-1\}}(-1)^{|I|}\Lambda_{(\mb{Y}_1,\hdots,\mb{Y}_k)_{I\cup\{k\}}}
\left(\frac{\boldsymbol{\tau}\,\mb{T}_1^{\boldsymbol{\beta}}}{\mb{r}} ,\hdots, \frac{\boldsymbol{\tau}\,\mb{T}_k^{\boldsymbol{\beta}}}{\mb{r}}
\right)_{I\cup\{k\}}\right).
\end{array}
\end{eqnarray}
\end{pro}
\dem Observe that
\begin{eqnarray}\nn
\begin{array}{rl}
&\dst n\, P\left(\mb{X}_1\leq \mb{u}_n^{(\boldsymbol{\tau})},\hdots,\mb{X}_{k-1}\leq \mb{u}_n^{(\boldsymbol{\tau})},\mb{X}_{k}\not\leq \mb{u}_n^{(\boldsymbol{\tau})}\right)\vspace{0.35cm}\\
=&\dst n \,\left\{P\left(\mb{X}_{k}\not\leq \mb{u}_n^{(\boldsymbol{\tau})}\right)-
\sum_{\emptyset\not=I\subset\{1,\hdots,k-1\}}(-1)^{|I|+1}
P\left(\bigcap_{i\in I}\mb{X}_1\not\leq \mb{u}_n^{(\boldsymbol{\tau})},\mb{X}_{k}\not\leq \mb{u}_n^{(\boldsymbol{\tau})}\right)\right\}\vspace{0.35cm}\\
=& \dst \sum_{\emptyset\subseteq I\subset\{1,\hdots,k-1\}}(-1)^{|I|}nP\left(\bigcap_{i\in I\cup\{k\}}\mb{Y}_i\not\leq { \mb{U}_{\mb{Y}_i}\left(\frac{\mb{r}\,n}{\boldsymbol{\tau}\,\mb{T}_i^{\boldsymbol{\beta}}} \right)} \right),
\end{array}
\end{eqnarray}
which leads to the result. \fdem\\


\begin{cor}
If $\mb{X}$ satisfies strong-mixing and  ${\tt D}^{(k)}\left(\mb{u}_n^{(\boldsymbol{\tau})}\right)$ conditions, for all $\boldsymbol{\tau}\in \er_+^d$, and
if (\ref{LambdaY1k}) holds for all $\mb{x}_1,\hdots,\mb{x}_k\in {\overline{\er}_+^{\,d}}$ then
\begin{eqnarray}\label{cmindextdk}
\theta_{\mb{X}}(\tau_1,\hdots,\tau_d)=1-\frac{\dst E\left(\sum_{\emptyset\not= I\subset\{1,\hdots,k-1\}}(-1)^{|I|+1}\Lambda_{(\mb{Y}_1,\hdots,\mb{Y}_k)_{I\cup\{k\}}}
\left(\frac{\boldsymbol{\tau}\,\mb{T}_1^{\boldsymbol{\beta}}}{\mb{r}} ,\hdots, \frac{\boldsymbol{\tau}\,\mb{T}_k^{\boldsymbol{\beta}}}{\mb{r}}
\right)_{I\cup\{k\}}\right)}{\dst E\left(\sum_{\emptyset\not= I\subset\{1,\hdots,d\}}(-1)^{|I|-1}\Lambda_{(\mb{Y}_1)_I}
\left(\frac{\boldsymbol{\tau}\,\mb{T}_1^{\boldsymbol{\beta}}}{\mb{r}}\right)_I\right)},
\end{eqnarray}

\begin{eqnarray}\label{cindextdk}
\theta_{X_j}=1-\frac{1}{r_j}E\left(\sum_{\emptyset\not= I\subset\{1,\hdots,k-1\}}(-1)^{|I|+1}\Lambda_{({Y}_{1j},\hdots,{Y}_{kj})_{I\cup\{k\}}}
\left(T_{1j}^{\beta_j},\hdots, T_{kj}^{\beta_j}\right)_{I\cup\{k\}}\right),
\end{eqnarray}
{ Moreover,  for all $\mb{x}\in {\er}_+^{\,d}$, we have
$$
\begin{array}{c}
\dst\lim_{n\to\infty}P\left(\mb{M}_n\leq \mb{U}_{\mb{X}}(n\mb{x})\right)
=\exp\left\{-\dst E\left(\sum_{\emptyset\subseteq I\subset\{1,\hdots,k-1\}}(-1)^{|I|}\Lambda_{(\mb{Y}_1,\hdots,\mb{Y}_k)_{I\cup\{k\}}}
\left(\frac{\,\mb{T}_1^{\boldsymbol{\beta}}}{\mb{r}\mb{x}} ,\hdots, \frac{\,\mb{T}_k^{\boldsymbol{\beta}}}{\mb{r}\mb{x}}
\right)_{I\cup\{k\}}\right)\right\}.
\end{array}\nn\\\nn
$$}
\end{cor}

The result points out that, however much is the sequencial dependence in the stationary sequence $\mb{T}$, the extremal index will be unit if $\mb{Y}_1,\hdots,\mb{Y}_k$ are tail independent. This characteristic is illustrated in Example \ref{emrarmax} with multivariate pRARMAX models, which are adjusted to a bivariate financial series in the last section.

\medskip

\begin{ex}\label{ex1}
Let $\{\mb{W}_n=(W_{n1},\hdots,W_{nd})\}_{n\geq 1}$ be an i.i.d.~sequence of random vectors with independent marginals  $W_{nj}$, $j=1\hdots,d$, having distribution function (d.f.) $F_{W_{nj}}(x)=(1-x^{-\beta_j})^{1/2}$, and $Y_{nj}=W_{n+1,j}\vee W_{nj}$, \jd. $\mb{Y}$ is $2$-dependent and thus we can apply the results on the calculation of the multivariate extremal index under the  condition ${\tt D}^{(2)}\left(\mb{u}_n^{(\boldsymbol{\tau})}\right)$. We have $\Lambda_{(\mb{Y}_1)_I}(\mb{x}_I)=0$ if $|I|>1$ and
\begin{eqnarray}\label{ex1lambda}
\Lambda_{(\mb{Y}_1,\mb{Y}_2)}(\mb{x}_1,\mb{x}_2)
=\dst\lim_{t\to\infty}tP\left(\mb{Y}_1\not\leq \mb{U}_{\mb{Y}_1}\left(\frac{t}{\mb{x}_1}\right),
\mb{Y}_2\not\leq \mb{U}_{\mb{Y}_2}\left(\frac{t}{\mb{x}_2}\right)\right)
=\sum_{j=1}^d\frac{1}{2}(x_{1j}\wedge x_{2j})\,.
\end{eqnarray}
Then, by (\ref{cmindextdk}),
\begin{eqnarray}\nn
\theta_{\mb{X}}(\tau_1,\hdots,\tau_d)=1-\dst\frac{\dst E\left(\Lambda_{(\mb{Y}_1,\mb{Y}_2)}
\left(\frac{\dst \boldsymbol{\tau}\,\mb{T}_1^{\boldsymbol{\beta}}}{\mb{r}} , \frac{\boldsymbol{\tau}\,\mb{T}_2^{\boldsymbol{\beta}}}{\mb{r}}
\right)\right)}{\dst\sum_{j=1}^dE\left(
\frac{{\tau_j}\,{T}_{1j}^{{\beta_j}}}{{r_j}}\right)}
=1-\dst\frac{\dst\frac{1}{2}\sum_{j=1}^dE\left(
\frac{{\tau_j}\,{T}_{1j}^{{\beta_j}}}{{r_j}}
\wedge \frac{{\tau_j}\,{T}_{2j}^{{\beta_j}}}{{r_j}}\right)}
{\dst\sum_{j=1}^d
{\tau_j}}
\end{eqnarray}
and
$$
\theta_j=1-\frac{1}{2r_j}E\left({T}_{1j}^{{\beta_j}}
\wedge {T}_{2j}^{{\beta_j}}\right)\,.
$$
Observe that the extremal indexes of the sequences  $\{Y_{nj}\}_{n\geq 1}$, \jd, are all equal to $1/2$ while the presence of the random factors may increase this value for different $\theta_j$, \jd. Therefore, each marginal sequence may have a different tendency for clustering of high values.
We also find that, since the marginals of $\mb{Y}_n$ are independent, the dependence structure of $\mb{T}_n$ does not affect the clustering of high values of  $\mb{X}_n$. As expected from (\ref{cmindextdk}) and illustrated in the examples of Section \ref{spc}, only some tail dependence of  $\mb{Y}_n$ allows to account the dependence of $T_{nj}$, \jd, on the value of $\theta_{\mb{X}}(\tau_1,\hdots,\tau_d)$.
\end{ex}

\begin{ex}\label{ex2}
Suppose that $\mb{Y}$ is $2$-dependent and $\{\mb{T}_{n}\}_{n\geq 1}$ is a sequence of independent vectors with independent marginals having Bernoulli distribution with mean $p$. By  (\ref{cmindextdk}), we have
\begin{eqnarray}\nn
\theta_{\mb{X}}(\tau_1,\hdots,\tau_d)=
1-\frac{\dst\sum_{\emptyset\not= I,J\subset\{1,\hdots,d\}}p^{|I|+|J|}(1-p)^{2d-|I|-|J|}
\Lambda_{((\mb{Y}_1)_I,(\mb{Y}_2)_{J})}
\left(\left(\frac{\boldsymbol{\tau}}{\mb{r}}\right)_I , \left(\frac{\boldsymbol{\tau}}{\mb{r}}
\right)_{J}\right)}{\dst\sum_{\emptyset\not= J\subset\{1,\hdots,d\}}p^{|J|}(1-p)^{d-|J|}
\left(-\log G_{\mb{Y}_J}\left(\left(\frac{\boldsymbol{\tau}}{\mb{r}}
\right)_{J}\right)\right)}.
\end{eqnarray}
In the particular case of the previous example in (\ref{ex1lambda}), the function above becomes
\begin{eqnarray}\nn
\theta_{\mb{X}}(\tau_1,\hdots,\tau_d)=
1-\frac{\dst\frac{p^2}{2}\sum_{j=1}^d
\frac{{\tau_j}}{{r_j}}}
{\dst\sum_{j=1}^d
{\tau_j}}
\end{eqnarray}
and $\theta_j=1-\frac{p^2}{2r_j}$, \jd.

In the Example \ref{ex4} ahead, we illustrate this choice of $\mb{T}_n$ to model lost values of $\mb{Y}_n$, considering that the marginals $Y_{nj}$, \jd, have total dependence.
\end{ex}

\section{Tail dependence}\label{sti}
The bivariate upper tail dependence of two random variables $X_i$ and $X_j$ can be measured through the tail dependence coefficient $\Lambda_{(X_i,X_j)}(1,1)$ (Sibuya \cite{sib} 1960, Joe \cite{joe} 1997). Extending this concept to cross-sectional lag-$m,s$ upper tail dependence, $m\geq 1$ and $0\leq s<d$, for a sequence $\mb{X}=\{\mb{X}_n\}_{n\geq 1}$, we have
\begin{eqnarray}\label{tdc-mr}
\begin{array}{c}
\lambda^{(m,s)}_X=\lim_{t\to\infty}\,t P\left(X_{1,1}>U_{X_{1}}(t),X_{1+m,1+s}>U_{X_{1+s}}(t)\right).
\end{array}
\end{eqnarray}
If $\lambda^{(m,s)}_X=0$ we say that $X_{1,1}$ and $X_{1+m,1+s}$ are upper tail independent. In this case, it is possible that a residual tail dependence captured at penultimate high levels may occur. This is measured through the asymptotic tail independent coefficient $\eta\in (0,1]$ (Ledford and Tawn \cite{led+tawn1,led+tawn2}, 1996/97). Analogously, we can extend this concept to cross-sectional lag-$m,s$:
\begin{eqnarray}\label{etamrxy}
\lim_{t\to\infty}\frac{P\left(X_{1,1}>U_{X_{1}}(t/x),X_{1+m,1+s}>U_{X_{1+s}}(t/y)\right)}
{P\left(X_{1,1}>U_{X_{1}}(t),X_{1+m,1+s}>U_{X_{1+s}}(t)\right)}=h_X^{(m,s)}(x,y),
\end{eqnarray}
for all $x,y\geq 0$, where $h_X^{(m,s)}$ is some non-degenerate function, homogeneous of order $-1/\eta_X^{(m,s)}$,  $\eta_X^{(m,s)}\in(0,1]$ and such that $h_X^{(m,s)}(1,1)=1$.
In particular, we have
\begin{eqnarray}\label{etamr}
P\left(X_{1,1}>U_{X_{1}}(t),X_{1+m,1+s}>U_{X_{1+s}}(t)\right)=t^{-1/\eta_X^{(m,s)}}l_X^{(m,s)}(t),
\end{eqnarray}
where $l_X^{(m,s)}$ is a slowly varying function, standing for the relative strength of dependence given a particular value of $\eta_X^{(m,s)}$. Whenever $\eta_X^{(m,s)}= 1 $ and $l_X^{(m,s)}$ converge to some constant $0<a\leq 1$ then  $X_{1,1}$ and $X_{1+m,1+s}$ are
tail dependent ($\lambda^{(m,s)}_X=a$), otherwise we have asymptotic tail independence with positive
association if $\eta_X^{(m,s)} > 1/2$, negative association if $\eta_X^{(m,s)} < 1/2$ and (almost) independence
if $\eta_X^{(m,s)} = 1/2$ (perfect if $l_X^{(m,s)}=1$).\\

By using the same arguments as in (\ref{utcopxy}), we obtain the following relation for $\lambda^{(m,s)}_X$.

\begin{pro}
The lag-$m,s$ upper tail dependence coefficient of $\mb{X}$ is given by
\begin{eqnarray}\label{tdc_cor}
\begin{array}{c}
\dst\lambda^{(m,s)}_X=
E\left(\Lambda_{(Y_{1,1},Y_{1+m,1+s})}\left(\frac{T_{1,1}^{\beta_1}}{r_1},\frac{T_{1+m,1+s}^{\beta_{1+s}}}{r_{1+s}}\right)\right),
\end{array}
\end{eqnarray}
where $\Lambda_{(Y_{1,1},Y_{1+m,1+s})} (x,y)=\lim_{t\to\infty}tP\left(Y_{1,1}>U_{Y_1}\left(t/x\right),Y_{1+m,1+s}>U_{Y_1+s}\left(t/y\right)\right)$ is the lag-$m,s$ upper tail copula of $\mb{Y}$.
\end{pro}

\begin{pro}\label{peta}
If $(Y_{1,1},Y_{1+m,1+s})$ satisfies (\ref{etamrxy}) for some function $h_Y^{(m,s)}$ and coefficient $\eta_Y^{(m,s)}$, then
$$
\begin{array}{c}
\dst P\left(X_{1,1}>U_{X_{1}}(t),X_{1+m,1+s}>U_{X_{1+s}}(t)\right)=
E\left(h_Y^{(m,s)}\left(\frac{T_{1,1}^{\beta_1}}{r_1}
,\frac{T_{1+m,1+s}^{\beta_{1+s}}}{r_{1+s}}\right)\right)
t^{-1/\eta_Y^{(m,s)}}l_Y^{(m,s)}(t)(1+o(1)),
\end{array}
$$
for large $t$.
\end{pro}
\dem Just observe that, using (\ref{etamrxy}), we have for large $t$
\begin{eqnarray}\nn
\begin{array}{rl}
&\dst P\left(X_{1,1}>U_{X_1}\left(t\right),X_{1+m,1+s}>U_{X_1+s}\left(t\right)\right)\vspace{0.35cm}\\
=& \dst \int P\left(Y_{1,1}>U_{Y_1}\left(\frac{r_1\,t}{z_{1}^{\beta_1}}\right),
Y_{1+m,1+s}>U_{Y_{1+s}}\left(\frac{r_{1+s}\,t}{z_{2}^{\beta_{1+s}}}\right)\right)dP_{(T_{1,1},T_{1+m,1+s})}(z_{1},z_{2})\vspace{0.35cm}\\
=&\dst P\left(Y_{1,1}>U_{Y_{1}}(t),Y_{1+m,1+s}>U_{Y_{1+s}}(t)\right)\int h_Y^{(m,s)}\left(\frac{z_{1}^{\beta_1}}{r_1},\frac{z_{2}^{\beta_{1+s}}}{r_{1+s}}\right)dP_{(T_{1,1},T_{1+m,1+s})}(z_{1},z_{2})(1+o(1)).
\end{array}
\end{eqnarray}
Now the result is straightforward from  (\ref{etamr}). \fdem\\

Observe that the lag-$m,s$ upper tail dependence or  independence of $\mb{X}$ is ruled by $\mb{Y}$. In the asymptotic tail independence case, this is even more evident since, by Proposition \ref{peta}, they have the same lag-$m,s$ asymptotic tail independent coefficient, i.e., $\eta_X^{(m,s)}=\eta_Y^{(m,s)}$, only differing in the relative strength  of dependence coming from the slowly varying functions.

If we consider the models in Examples \ref{ex1} and \ref{ex2}, we derive $\lambda^{(m,s)}_X=0$, for all $m\geq 1$ and $1\leq s <d$. Observe that this corresponds to perfect independence of $\mb{Y}$ and thus it satisfies (\ref{etamrxy}) and (\ref{etamr}), with $\eta_Y^{(m,s)}=1/2$, $h_Y^{(m,s)}(x,y)=xy$ and $l_Y^{(m,s)}(t)=1$. Therefore,
\begin{eqnarray}\label{exetamrxy}
\begin{array}{c}
P\left(X_{1,1}>U_{X_{1}}(t),X_{1+m,1+s}>U_{X_{1+s}}(t)\right)=t^{-2},
\end{array}
\end{eqnarray}
i.e., perfect independence of $\mb{X}$ too.


\section{Some particular cases that extend existing models}\label{spc}

In this section we consider our model restricted to a common rescaled factor $Y$, i.e., by taking $Y_{nj}=Y_n^{\gamma_j}$, $\gamma_j>0$, \jd, with $\mb{Y}^{(1)}=\{Y_{n}\}_{n\geq 1}$ a stationary sequence having common Pareto-type d.f. with shape parameter $\alpha>0$, and $T_{nj}=Z_{nj}^{\gamma_j}$, \jd, where $\mb{Z}=\{(Z_{n1}, \hdots, Z_{nd})\}_{n\geq 1}$ is a stationary sequence, independent of $\mb{Y}^{(1)}$, with support $\er_+^d$ and such that, $r_j=E(Z_{nj}^{\alpha})<\infty$, $j=1,\hdots, d$. We will therefore apply the previous results with $\beta_j=\alpha/\gamma_j$, $j=1,\hdots,d$ and $\Lambda_{\mb Y}(x_1,\hdots,x_d)=\bigwedge_{j=1}^dx_j$.\\

This case extends the heavy-tailed factor model in Lescourret and Robert (\cite{les+rob06}, 2006), which was used to measure the extreme dependence in loss severities of storm insurance data and where $Y$ was considered a common latent factor corresponding to the intensity of the natural disaster. In addition, if we consider $\gamma=\gamma_j=1$, $j=1,\hdots,d$, we obtain an extended version of the scale mixture model of Li (\cite{li}, 2009) with several applications to real data (see, e.g., Arnold \cite{arn83} 1983, Kotz \emph{et al.} \cite{kotz+2000} 2000 and references therein).\\

As a consequence of Proposition \ref{putc}, we obtain the following result for the upper tail copula of $\mb{X}$:
\begin{eqnarray}\nn
\Lambda_{\mb{X}}(x_1,\hdots,x_d)
=E\left(\bigwedge_{j=1}^d\frac{x_jT_{j}^{\beta_j}}{r_j}\right)
=E\left(\bigwedge_{j=1}^d\frac{x_jZ_{j}^{\alpha}}{r_j}\right),
\end{eqnarray}
for all $(x_1,...,x_d)\in\overline{\er}_+^d$, where $(Z_1,\hdots,Z_d)\dst\mathop{=}^d(Z_{n1},\hdots,Z_{nd})$. From Proposition \ref{pmev}, we conclude that $F_{\mb{X}_1}$ is in the domain of attraction of
\begin{eqnarray}\label{CRMEV}
G(x_1,\hdots,x_d)=\exp\left\{-E\left(\bigvee_{j=1}^d\frac{x_j^{-1}Z_{nj}^{\alpha}}{r_j}\right)\right\}.
\end{eqnarray}

\begin{rem}
If we assume that $\mb{Z}$ is an i.i.d.\,sequence of random vectors with i.i.d.\,marginals $Z_{nj}$, $j=1\hdots,d$, having common d.f.~$H$ and $r_j=r$, \jd, we easily derive
\begin{eqnarray}\nn
\begin{array}{c}
E\left(\bigwedge_{j=1}^d\frac{x_j\,{Z}_{1j}^{{\alpha}}}{r_j} \right)= \frac{1}{r}\sum_{j=1}^dx_j \, E\left(\frac{Z^\alpha}{1-H(Z)}\prod_{i=1}^d\left[1- H\left((x_j/x_i)^{1/\alpha}\,Z\right)\right]\right)
\end{array}
\end{eqnarray}
and
\begin{eqnarray}\nn
\begin{array}{c}
E\left(\bigvee_{j=1}^d\frac{x_j^{-1}\,{Z}_{1j}^{{\alpha}}}{r_j} \right)= \frac{1}{r}\sum_{j=1}^dx_j^{-1} \, E\left(\frac{Z^\alpha}{H(Z)}\prod_{i=1}^d H\left((x_i/x_j)^{1/\alpha}\,Z\right)\right),
\end{array}
\end{eqnarray}
extending the bivariate results in Lescourret and Robert (\cite{les+rob06}, 2006).\\
\end{rem}

In the sequel, we rewrite Proposition \ref{pdk}, compute the extremal index and find the attractor MEV within this particular model.

In what concerns the local dependence condition of Proposition \ref{pdk}, observe that
$$
\left\{\mb{Y}^{(1)}_i\not\leq\mb{v}_n^{(\boldsymbol{\tau}^*(\mb{z}))}\right\}
=\left\{Y_i>\bigwedge_{j=1}^d\left(\frac{r_jn}{z_{ij}^\alpha\tau_j}\right)^{1/\alpha}\right\},
$$
$$
\left\{\mb{Y}^{(1)}_i\leq\mb{v}_n^{(\boldsymbol{\tau}^*(\mb{z}))}\right\}
=\left\{Y_i\leq \bigwedge_{j=1}^d\left(\frac{r_jn}{z_{ij}^\alpha\tau_j}\right)^{1/\alpha}\right\}
$$
and
$$
a^\alpha\bigvee_{j=1}^d\frac{\tau_j}{r_j}
<\bigvee_{j=1}^d\frac{z_{ij}^\alpha\tau_j}{r_j}
<b^\alpha\bigvee_{j=1}^d\frac{\tau_j}{r_j}.
$$
Therefore, we obtain for this model the following particular result.
\begin{pro}\label{pCRdk}
If $Z_{nj}$ has support in $[a,b]$, $j=1,\hdots,d$, and  $\mb{Y}^{(1)}$ satisfies condition ${\tt D}^{(k)}\left({v}_n^{({\tau}({a}))},{v}_n^{({\tau}({b}))}\right)$
defined by
$$
\lim_{n\to\infty} n\sum_{i=k}^{[n/k_n]}P\left({Y}_{1}> {v}_n^{({\tau}({b}))},{Y}_{i}\leq {v}_n^{({\tau}({a}))},{Y}_{i+1}> {v}_n^{({\tau}({b}))}\right)=0,
$$
with ${\tau}({a})=\bigvee_{j=1}^d\frac{a^\alpha\tau_j}{r_j}>
{\tau}({b})=\bigvee_{j=1}^d\frac{b^\alpha\tau_j}{r_j}$, $\tau_j>0$,  $j=1,\hdots,d$, then ${\tt D}^{(k)}\left(\mb{u}_n^{(\boldsymbol{\tau})}\right)$ holds for $\mb{X}$, for each $\boldsymbol{\tau}=(\tau_1,\hdots,\tau_d)$.\\
\end{pro}

Assuming that
\begin{eqnarray}\label{LambdaY1kuniv}
\Lambda_{({Y}_1,\hdots,{Y}_k)}({x}_1,\hdots,{x}_k)=\dst\lim_{t\to\infty}
t\,P\left({Y}_i>\bigcap_{i=1}^k\left\{ U_{Y_i}\left(t/x_i\right)\right\} \right),
\end{eqnarray}
exists with $(x_1,\hdots,x_k)\in\overline{\er}_+^{\,k}$, then
\begin{eqnarray}\nn
\begin{array}{rl}
&\dst\lim_{n\to\infty} n\, P\left(\mb{X}_1\leq \mb{u}_n^{(\boldsymbol{\tau})},\hdots,\mb{X}_{k-1}\leq \mb{u}_n^{(\boldsymbol{\tau})},\mb{X}_{k}\not\leq \mb{u}_n^{(\boldsymbol{\tau})}\right)\vspace{0.35cm}\\
=&\dst E\left(\sum_{\emptyset\subseteq I\subset\{1,\hdots,k-1\}}(-1)^{|I|}\Lambda_{({Y}_1,\hdots,{Y}_k)_{I\cup\{k\}}}
\left(\bigvee_{j=1}^d\frac{{\tau_j}\,{Z}_{1j}^{{\alpha}}}{{r_j}} ,\hdots, \bigvee_{j=1}^d\frac{{\tau_j}\,{Z}_{kj}^{{\alpha}}}{{r_j}}
\right)_{I\cup\{k\}}\right),
\end{array}
\end{eqnarray}
with the convention that $\Lambda_{{Y}_k}\left(\bigvee_{j=1}^d\frac{{\tau_j}\,{Z}_{kj}^{{\alpha}}}{{r_j}}\right)
=\bigvee_{j=1}^d\frac{{\tau_j}\,{Z}_{kj}^{{\alpha}}}{{r_j}}$. Therefore, if $\mb{X}$ satisfies strong-mixing and  ${\tt D}^{(k)}\left(\mb{u}_n^{(\boldsymbol{\tau})}\right)$ conditions, for all $\boldsymbol{\tau}\in \er_+^d$, and
if (\ref{LambdaY1kuniv}) holds for all $(x_1,\hdots,x_k)\in\overline{\er}_+^{\,k}$, then
\begin{eqnarray}\label{CRmindextdk}
\theta_{\mb{X}}(\tau_1,\hdots,\tau_d)=1-\dst\frac{\dst E\left(\sum_{\emptyset\not= I\subset\{1,\hdots,k-1\}}(-1)^{|I|+1}\Lambda_{({Y}_1,\hdots,{Y}_k)_{I\cup\{k\}}}
\left(\bigvee_{j=1}^d\frac{{\tau_j}\,{Z}_{1j}^{{\alpha}}}{{r_j}} ,\hdots, \bigvee_{j=1}^d\frac{{\tau_j}\,{Z}_{kj}^{{\alpha}}}{{r_j}}\right)_{I\cup\{k\}}\right)}
{\dst E\left(\bigvee_{j=1}^d\frac{{\tau_j}\,{Z}_{1j}^{{\alpha}}}{{r_j}} \right)},
\end{eqnarray}
the marginal extremal index is, for $j=1,\hdots,d$,
\begin{eqnarray}\label{CRindextdk}
\theta_{X_j}=1-\frac{1}{r_j}\, E\left(\sum_{\emptyset\not= I\subset\{1,\hdots,k-1\}}(-1)^{|I|+1}\Lambda_{({Y}_1,\hdots,{Y}_k)_{I\cup\{k\}}}
\left(Z_{1j}^{\alpha} ,\hdots, Z_{kj}^{\alpha}\right)_{I\cup\{k\}}\right)
\end{eqnarray}
and, as $n\to\infty$, $P\left(\mb{M}_n\leq \mb{U}_{\mb{X}}(n\mb{x})\right)$ converges to the limiting MEV
\begin{eqnarray}\nn
\exp\left\{-E\left(\sum_{\emptyset\subseteq I\subset\{1,\hdots,k-1\}}(-1)^{|I|}\Lambda_{({Y}_1,\hdots,{Y}_k)_{I\cup\{k\}}}
\left(\bigvee_{j=1}^d\frac{{Z}_{1j}^{{\alpha}}}{x_j\,r_j} ,\hdots, \bigvee_{j=1}^d\frac{{Z}_{kj}^{{\alpha}}}{x_j\,r_j}
\right)_{I\cup\{k\}}\right)\right\},
\end{eqnarray}
for all $(x_1,\hdots,x_d)\in{\er}_+^{\,d}$.\\


The cross-sectional lag-$m,s$ upper tail dependent coefficient, $m\geq 1$ and $0\leq s<d$, is given by
\begin{eqnarray}\label{CRtdc_cor}
\begin{array}{c}
\dst \lambda^{(m,s)}_X=
E\left(\Lambda_{(Y_{1},Y_{1+m})}\left(\frac{Z_{1,1}^{\alpha}}{r_1}, \frac{Z_{1+m,1+s}^{\alpha}}{r_{1+s}}\right)\right),
\end{array}
\end{eqnarray}
In case $\mb{X}$ is cross-sectional lag-$m,s$ upper tail independent (i.e., $\lambda^{(m,s)}_X=0$) and if $(Y_{1},Y_{1+m})$ satisfies (\ref{etamrxy}) for some function $h_Y^{(m)}$ and asymptotic tail independent coefficient $\eta_Y^{(m)}$, for large $t$,  then
\begin{eqnarray}\label{CRetamrxy}
\begin{array}{c}
\dst P\left(X_{1,1}>U_{X_{1}}(t),X_{1+m,1+s}>U_{X_{1+s}}(t)\right)=
E\left(h_Y^{(m)}\left(\frac{Z_{1,1}^{\alpha}}{r_1}
,\frac{Z_{1+m,1+s}^{\alpha}}{r_{1+s}}\right)\right)
t^{-1/\eta_Y^{(m)}}l_Y^{(m)}(t)(1+o(1)).
\end{array}
\end{eqnarray}

Next, we consider some particular examples.

\begin{ex}\label{emrarmax}[Multivariate autoregressive  processes with random coefficients]
Suppose that $Y_n=c(Y_{n-1}\vee W_n)$, $n\geq 1$, where $0<c<1$ is a constant and $\{W_n\}_{n\geq 1}$ is an i.i.d. sequence, independent of $Y_0$, with common d.f. $F_W$. The max-autoregressive sequence $\mb{Y}^{(1)}= \{Y_{n}\}_{n\geq 1} $, usually denoted ARMAX, has stationary distribution $F_Y$ if and only if $0<\sum_{j=1}^\infty(1-F_W(x/c^j))<\infty$, for some $x>0$ and, in this case, $F_Y(x)=F_Y(x/c)F_W(x/c)$ (Alpuim, \cite{alp1} 1989). The Pareto distribution with parameter $\alpha>0$ is a stationary d.f. of $\mb{Y}$ with $F_W(x)=(1-(cx)^{-\alpha})/(1-x^{-\alpha})$.

If $Z_{nj}$ has support in $[a,b]$, $j=1,\hdots,d$, such that $c<a/b $, then $\mb{Y}^{(1)}$ satisfies condition ${\tt D}^{(2)}\left({v}_n^{({\tau}({a}))},{v}_n^{({\tau}({b}))}\right)$, with $\tau(a)$ and  $\tau(b)$ as in Proposition \ref{pCRdk}. Indeed, we have
\begin{eqnarray}\nn
\begin{array}{rl}
&\dst n\,\sum_{i=2}^{[n/k_n]}P(Y_1>v_n^{(\tau(b))},Y_i\leq  v_n^{(\tau(a))}, Y_{i+1}>  v_n^{(\tau(b))})\vspace{0.35cm}\\
\leq&\dst n\,\sum_{i=2}^{[n/k_n]}P(Y_1>v_n^{(\tau(b))},Y_i\leq  v_n^{(\tau(a))}, Y_{i}>  v_n^{(\tau(b))}/c)\vspace{0.15cm}\\
&\dst +n\,\sum_{i=2}^{[n/k_n]}P(Y_1>v_n^{(\tau(b))},Y_i\leq  v_n^{(\tau(a))}, W_{i+1}>  v_n^{(\tau(b))}/c).
\end{array}
\end{eqnarray}
The first term is null since $c<a/b\Leftrightarrow v_n^{(\tau(a))}<u_n^{(\tau(b))}/c$ and the second term is upper-bounded, successively, by
\begin{eqnarray}\nn
\begin{array}{rl}
&\dst n\,\frac{n}{k_n}P(Y_1>v_n^{(\tau(b))}P(W_{i+1}>  v_n^{(\tau(b))})= n\,\frac{n}{k_n}(1-F_{Y}(v_n^{(\tau(b))}))
\left(1-\frac{F_{Y_1}(v_n^{(\tau(b))})}{F_{Y_1}(v_n^{(\tau(b))}/c)}\right)\vspace{0.35cm}\\
\leq&\dst \frac{n^2}{k_n}(1-F_{Y}(v_n^{(\tau(b))}))^2\to 0
\end{array}
\end{eqnarray}
for any $k_n\to\infty$, as $n\to\infty$. Therefore, if the ARMAX sequence ${\mb Y}$ is such that $c<a/b$ and ${\mb Z}$ has support in $[a,b]^d$, then ${\mb X}$ satisfies ${\tt D}^{(2)}\left(\mb{u}_n^{(\boldsymbol{\tau})}\right)$, for all $\boldsymbol{\tau}\in \er_+^d$ (Proposition \ref{pCRdk}).


From Ferreira and Canto e Castro (\cite{mf+lcc}, 2008; Proposition 3.6), we obtain the following bivariate upper tail copula function,  for $(x,y)\in\overline{\er}_+^d$:
\begin{eqnarray}\nn
\Lambda_{({Y}_1,Y_{1+m})}(x,y)=x\wedge yc^{m\alpha}
\end{eqnarray}

Now we assume that $\mb{Z}$ is an i.i.d.\,sequence of random vectors with i.i.d.\,marginals $Z_{nj}$, $j=1\hdots,d$, having common d.f.\,$H$ and $r_j=r$, \jd. After some simple calculations, we derive
\begin{eqnarray}\nn
\begin{array}{rl}
&\dst E\left(\bigvee_{j=1}^d\frac{x_j\,{Z}_{1j}^{{\alpha}}}{r_j} \wedge c^{\alpha}\bigvee_{j=1}^d\frac{x_j\,{Z}_{2j}^{{\alpha}}}{r_j}\right)\vspace{0.35cm}\\
=&\dst  \frac{1}{r}\sum_{j=1}^dx_j \, E\left(\frac{Z^\alpha}{H(Z)}\prod_{i=1}^d H\left((x_j/x_i)^{1/\alpha}\,Z\right)\left[1-\prod_{i=1}^d H\left((x_j/x_i)^{1/\alpha}\,Z/c\right)\right]\right)\vspace{0.35cm}\\
&\dst +\frac{1}{r}\sum_{j=1}^d x_j \, E\left(\frac{(cZ)^\alpha}{H(Z)}\prod_{i=1}^d H\left((x_j/x_i)^{1/\alpha}\,Z\right)\left[1-\prod_{i=1}^d H\left((x_j/x_i)^{1/\alpha}\,Zc\right)\right]\right),
\end{array}
\end{eqnarray}
where $r=r_j$, \jd.

Therefore, $\mb{X}$ has multivariate extremal index given by
\begin{eqnarray}\nn
\begin{array}{rl}
\dst \theta_{\mb{X}}(\tau_1,\hdots,\tau_d)=1-\frac{\sum_{j=1}^d\tau_j \left\{ \, E\left(\frac{Z^\alpha}{H(Z)}\pi(Z)\left[1-\pi(Z/c)\right]\right)
+E\left(\frac{(cZ)^\alpha}{H(Z)}\pi(Z)\left[1-\pi(Zc)\right]\right)\right\}}
{\sum_{j=1}^d\tau_j \, E\left(\frac{Z^\alpha}{H(Z)}\pi(Z)\right)}
\end{array}
\end{eqnarray}
where $\pi(y)=\prod_{i=1}^d H\left((\tau_j/\tau_i)^{1/\alpha}\,y\right)$, and for each \jd,
\begin{eqnarray}\nn
\begin{array}{rl}
\dst \theta_{{X_j}}=1-\frac{ E\left(Z^\alpha\left(1-H(Z/c)\right)
+(cZ)^\alpha\left(1-H(Zc)\right)\right)}
{ r}.
\end{array}
\end{eqnarray}

In what concerns the limiting MEV of $\mb{X}$, for all $(x_1,\hdots,x_d)\in{\er}_+^{\,d}$, we have
\begin{eqnarray}\nn
\begin{array}{c}
\dst\lim_{n\to\infty}P\left(\mb{M}_n\leq \mb{U}_{\mb{X}}(n\mb{x})\right)=\exp\left\{-\frac{1}{r}\sum_{j=1}^dx_j^{-1} \left[ E\left(\frac{Z^\alpha}{H(Z)}\pi^*(Z)\pi^*(Z/c)\right) -E\left(\frac{(cZ)^\alpha}{H(Z)}\pi^*(Z)\left[1-\pi^*(Zc)\right]\right)\right]\right\}.
\end{array}
\end{eqnarray}
where $\pi^*(y)=\prod_{i=1}^d H\left((x_i/x_j)^{1/\alpha}\,y\right)$.

Analogously, we obtain the cross-sectional lag-$m,s$  upper tail dependent coefficient given by
\begin{eqnarray}\nn
\begin{array}{c}
\dst \lambda^{(m,s)}_X=\frac{1}{r}\,E\left(Z^\alpha H(c^mZ)
+(c^mZ)^\alpha\left(1-H(c^mZ)\right)\right).
\end{array}
\end{eqnarray}
Particular cases can be obtained by replacing $H$ by any d.f.\,with support in $[a,b]$, such that $c<a/b$.\\

Now if we consider, $Y_n=Y_{n-1}^c\vee W_n$, $n\geq 1$, with $0<c<1$ and $\{W_n\}_{n\geq 1}$ an i.i.d. sequence with common d.f. $F_W$, independent of $Y_0$, we obtain a max-autoregressive sequence with an exponent transformation, denoted pARMAX (Ferreira and Canto e Castro, \cite{mf+lcc}, 2008). A stationary distribution $F_Y$ exists if and only if $0<\sum_{j=1}^\infty(1-F_W(x^{c^j}))<\infty$, for some $x>0$,  and is such that $F_Y(x)=F_Y(x^{1/c})F_W(x)$. This relation is satisfied by a Pareto($\alpha$) distribution with $F_W(x)=(1-x^{-\alpha})/(1-x^{-\alpha/c})$ (see Ferreira and Canto e Castro \cite{lccmf2} 2010 and references therein). The pARMAX sequence $\mb{Y}^{(1)}$ has lag-$m$ upper tail dependence function $\Lambda_{(Y_1,Y_{1+m})}(x,y)=0$, for all $m\geq 1$, and thus it is an upper tail independent process. Therefore, in this case, the lag-$m,s$ upper tail dependence function $\Lambda_{(X_{1,1},X_{1+m,1+s})}(x,y)$ is also null, for all  $m\geq 1$ and $s=1,\hdots,d$, the extremal index is unit and the limiting MEV is as in the i.i.d.\,case, i.e., it is given by (\ref{CRMEV}).
Moreover, since $\mb{Y}^{(1)}$ satisfies (\ref{etamrxy}) with asymptotic tail independent coefficient $\eta_Y^{(m)}=1/2\vee c^m$ and $h_Y^{(m)}(x,y)=xy\mathds{1}_{\{c^m\leq 1/2\}}+y^{1/c^m}\mathds{1}_{\{c^m> 1/2\}}$, $x,y \geq 0$ (Ferreira and Canto e Castro, \cite{mf+lcc}, 2008), by applying (\ref{CRetamrxy}), we obtain \begin{eqnarray}\nn
\begin{array}{rl}
&\dst P\left(X_{1,1}>U_{X_{1}}(t),X_{1+m,1+s}>U_{X_{1+s}}(t)\right)\vspace{0.35cm}\\
=&\dst
\left(\mathds{1}_{\{c^m\leq 1/2\}}+
\frac{E\left(Z_{1+m,1+s}^{\alpha/c^m}\right)}{r_{1+s}^{1/c^m}}
\mathds{1}_{\{c^m> 1/2\}}\right)
t^{-1/(1/2\vee c^m)}l_Y^{(m)}(t)(1+o(1)).
\end{array}
\end{eqnarray}
Observe that sequence $\mb{X}$ generated from an ARMAX recursion $\mb{Y}^{(1)}$ corresponds to a multivariate formulation of the RARMAX process introduced in Alpuim and Athayde (\cite{alp+atha90}, 1990), with applications within  reliability and various natural phenomena. If sequence $\mb{X}$ is generated from a pARMAX recursion $\mb{Y}$, we have a multivariate formulation of the pRARMAX process introduced in Ferreira and Canto e Castro (\cite{lccmf2}, 2010), used in the modeling of financial series.
\end{ex}

\begin{ex}[Multivariate moving maxima processes with random coefficients]\label{ex4}
Consider $Y_n=W_{n+1}\vee W_n)$, $n\geq 1$, where $W=\{W_n\}_{n\geq 1}$ is an i.i.d. sequence with common d.f. $F_W(x)=(1-x^{-\alpha})^{1/2}$, $x\geq 1$. The sequence $\mb{Y}^{(1)}$ is $2$-dependent and thus $\mb{X}$ satisfies condition ${\tt D}^{(2)}(\mb{u}_n^{\boldsymbol{\tau}})$, for all $\boldsymbol{\tau}\in \er_+^d$, as well as, $\Lambda_{(Y_1,Y_{1+m})}(x,y)=0$, for $m\geq 2$. It is easily seen that
$$
\Lambda_{(Y_1,Y_{2})}(x,y)=\frac{1}{2}\left(x+y-x\vee y\right)=\frac{1}{2}\left(x\wedge y\right)
$$

Analogously to the ARMAX example above, we assume that $\mb{Z}$ is an i.i.d.\,sequence of random vectors with i.i.d.\,marginals $Z_{nj}$, $j=1\hdots,d$, having common d.f.\,$H$ and $r_j=r$, \jd, and a similar  procedure lead us to
\begin{eqnarray}\nn
\begin{array}{rl}
&\dst E\left(\frac{1}{2}\left(\bigvee_{j=1}^d\frac{x_j\,{Z}_{1j}^{{\alpha}}}{r_j} \wedge \bigvee_{j=1}^d\frac{x_j\,{Z}_{2j}^{{\alpha}}}{r_j}\right)\right)\vspace{0.35cm}\\
=& \dst \frac{1}{r}\sum_{j=1}^dx_j \, E\left(\frac{Z^\alpha}{H(Z)}\prod_{i=1}^d H\left((x_j/x_i)^{1/\alpha}\,Z\right)\left[1-\prod_{i=1}^d H\left((x_j/x_i)^{1/\alpha}\,Z\right)\right]\right),
\end{array}
\end{eqnarray}
where $r=r_j$, \jd. Hence, we derive the multivariate extremal index
\begin{eqnarray}\nn
\begin{array}{rl}
\dst \theta_{\mb{X}}(\tau_1,\hdots,\tau_d)=1-\frac{\sum_{j=1}^d\tau_j \, E\left(\frac{Z^\alpha}{H(Z)}\prod_{i=1}^d H\left((\tau_j/\tau_i)^{1/\alpha}\,Z\right)\left[1-\prod_{i=1}^d H\left((\tau_j/\tau_i)^{1/\alpha}\,Z\right)\right]\right)}
{\sum_{j=1}^d\tau_j \, E\left(\frac{Z^\alpha}{H(Z)}\prod_{i=1}^d H\left((\tau_j/\tau_i)^{1/\alpha}\,Z\right)\right)},
\end{array}
\end{eqnarray}
the univariate extremal index, for \jd,
\begin{eqnarray}\nn
\begin{array}{rl}
\dst \theta_{{X_j}}=1-\frac{ E\left(Z^\alpha [1-H(Z)]\right)}
{ r},
\end{array}
\end{eqnarray}
and, for the limiting MEV, with $(x_1,\hdots,x_d)\in{\er}_+^{\,d}$,
\begin{eqnarray}\nn
\begin{array}{c}
\dst\lim_{n\to\infty}P\left(\mb{M}_n\leq \mb{U}_{\mb{X}}(n\mb{x})\right) =\exp\left\{-\frac{1}{r}\sum_{j=1}^dx_j^{-1} E\left(\frac{Z^\alpha}{H(Z)}\left[\prod_{i=1}^d H\left((x_i/x_j)^{1/\alpha}\,Z\right)\right]^2\right)
\right\}.
\end{array}
\end{eqnarray}

The cross-sectional lag-$1,s$  upper tail dependent coefficient is
\begin{eqnarray}\nn
\begin{array}{c}
\dst \lambda^{(1,s)}_X=\frac{1}{2}=\lambda^{(1)}_Y.
\end{array}
\end{eqnarray}
The case $m\geq 2$ corresponding to perfect independence of $\mb{Y}^{(1)}$ has been analyzed in (\ref{exetamrxy}).\\


Observe that each marginal $X_{n,j}$ is a moving maxima process with a random coefficient $Z_{n,j}$, \jd. Therefore $\mb{X}$ corresponds to a multivariate formulation of this latter. For applications of multivariate moving maxima processes, see Zhang (\cite{zhang09}, 2009) and references therein.\\

We are going to consider some particular cases of  $H$.\\
\begin{itemize}
\item({Multivariate Gumbel case}) Suppose that $H$ is a Fr\'echet($\delta$, $\xi$) distribution, $\delta,\xi>0$, i.e., $H(x)=\exp(-\delta\,x^{-\xi})$, $x>0$, with $\xi>\alpha$. Therefore, $P(Z^\alpha\leq x)=\exp(-\delta\,x^{-\xi/\alpha})$ and $E(Z^\alpha)=e^{\alpha/\xi}\Gamma(1-\alpha/\xi)<\infty$, since $\xi/\alpha>1$. We have
\begin{eqnarray}\nn
\begin{array}{c}
\dst E\left(\frac{Z^\alpha}{H(Z)}\prod_{i=1}^d H\left((x_i/x_j)^{1/\alpha}\,Z\right)\right)=
\frac{1}{\dst \sum_{i=1}^d(x_j/x_i)^{\xi/\alpha}}\left[\delta\sum_{i=1}^d(x_j/x_i)^{\xi/\alpha}\right]^{\alpha/\xi}
\Gamma(1-\alpha/\xi),
\end{array}
\end{eqnarray}
as well as
\begin{eqnarray}\nn
\begin{array}{c}
\dst E\left(\frac{Z^\alpha}{H(Z)}\left[\dst \prod_{i=1}^d H\left((x_i/x_j)^{1/\alpha}\,Z\right)\right]^2\right)=
\frac{1}{\dst 2\sum_{i=1}^d(x_j/x_i)^{\xi/\alpha}}\left[2\delta\sum_{i=1}^d(x_j/x_i)^{\xi/\alpha}\right]^{\alpha/\xi}
\Gamma(1-\alpha/\xi).
\end{array}
\end{eqnarray}
Therefore,
\begin{eqnarray}\nn
\begin{array}{c}
\dst G(x_1,\hdots,x_d)=\exp\left\{-\left(\sum_{i=1}^dx_i^{-\xi/\alpha}\right)^{\alpha/\xi}\right\},
\end{array}
\end{eqnarray}
which is a multivariate Gumbel distribution, and
\begin{eqnarray}\nn
\begin{array}{c}
\dst\lim_{n\to\infty}P\left(\mb{M}_n\leq \mb{U}_{\mb{X}}(n\mb{x})\right) =\exp\left\{-2^{\alpha/\xi-1}\left(\sum_{i=1}^dx_i^{-\xi/\alpha}\right)^{\alpha/\xi}
\right\}.
\end{array}
\end{eqnarray}
Thus the extremal index is constant and is given by
$$
\theta_{\mb X}(\tau_1,\hdots,\tau_d)=2^{\alpha/\xi-1}=\theta_j,\,j=1\hdots,d.
$$
\item({Thinning model}) Suppose that $H$ is a Bernoulli($p$) distribution, $0<p<1$. Since
\begin{eqnarray}\nn
\begin{array}{c}
\dst E\left(\frac{Z^\alpha}{H(Z)}\prod_{i=1}^d H\left((x_i/x_j)^{1/\alpha}\,Z\right)\right)=
p(1-p)^{\sum_{i=1}^d\mathds{1}_{\{x_i<x_j\}}}
\end{array}
\end{eqnarray}
and
\begin{eqnarray}\nn
\begin{array}{c}
\dst E\left(\frac{Z^\alpha}{H(Z)}\left[\prod_{i=1}^d H\left((x_i/x_j)^{1/\alpha}\,Z\right)\right]^2\right)=
p(1-p)^{2\sum_{i=1}^d\mathds{1}_{\{x_i<x_j\}}},
\end{array}
\end{eqnarray}
we have
\begin{eqnarray}\nn
\begin{array}{c}
\dst G(x_1,\hdots,x_d)=\exp\left\{-\sum_{j=1}^dx_j^{-1}
(1-p)^{\sum_{i=1}^d\mathds{1}_{\{x_i<x_j\}}}\right\},
\end{array}
\end{eqnarray}
as well as
\begin{eqnarray}\nn
\begin{array}{c}
\dst\lim_{n\to\infty}P\left(\mb{M}_n\leq \mb{U}_{\mb{X}}(n\mb{x})\right) =\exp\left\{-\sum_{j=1}^dx_j^{-1}
(1-p)^{2\sum_{i=1}^d\mathds{1}_{\{x_i<x_j\}}}\right\}.
\end{array}
\end{eqnarray}
The extremal index is given by the function
$$
\dst \theta_{\mb X}(\tau_1,\hdots,\tau_d)=1-\frac{\dst \sum_{j=1}^d\tau_j
(1-p)^{2\sum_{i=1}^d\mathds{1}_{\{\tau_i>\tau_j\}}}}{\dst \sum_{j=1}^d\tau_j
(1-p)^{\sum_{i=1}^d\mathds{1}_{\{\tau_i>\tau_j\}}}}
$$
and $\theta_j=1$, \jd.
\end{itemize}

\end{ex}

\section{An application to financial data}\label{saplic}

\begin{figure}[!thb]
\centering
\includegraphics[width= 5.8cm,height= 3.9cm]{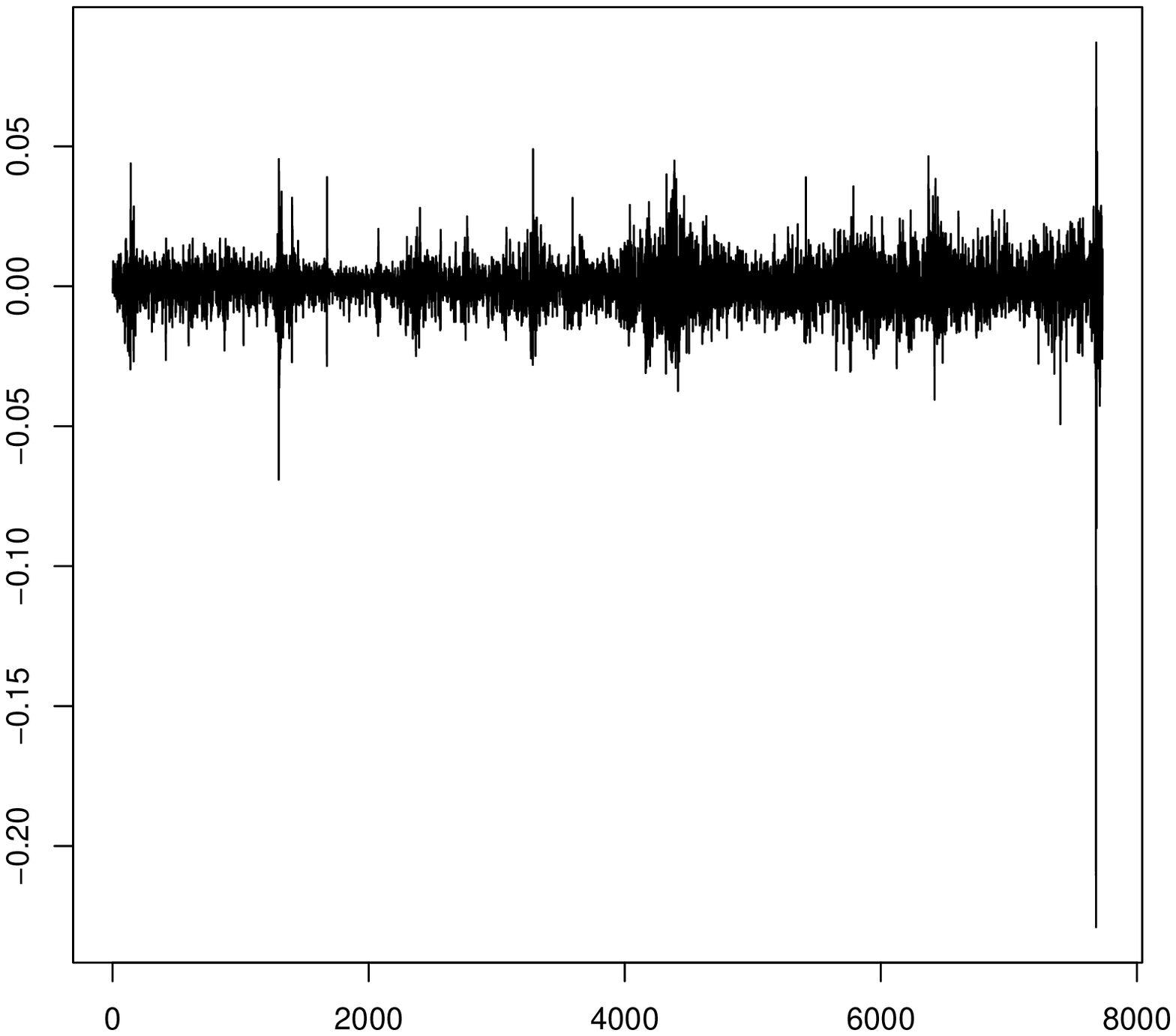}
\includegraphics[width= 5.8cm,height= 3.9cm]{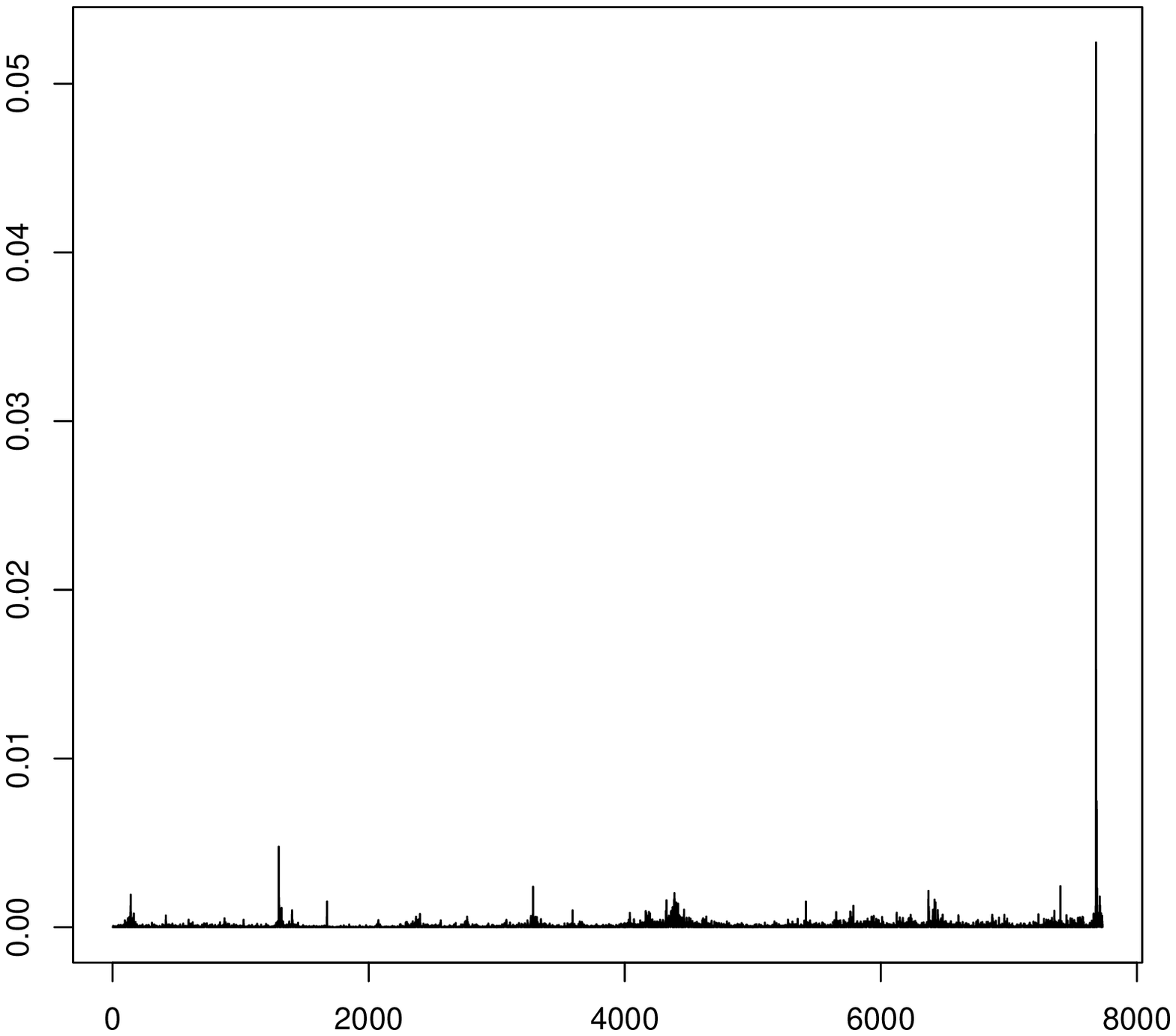}\\
\includegraphics[width= 5.8cm,height= 3.9cm]{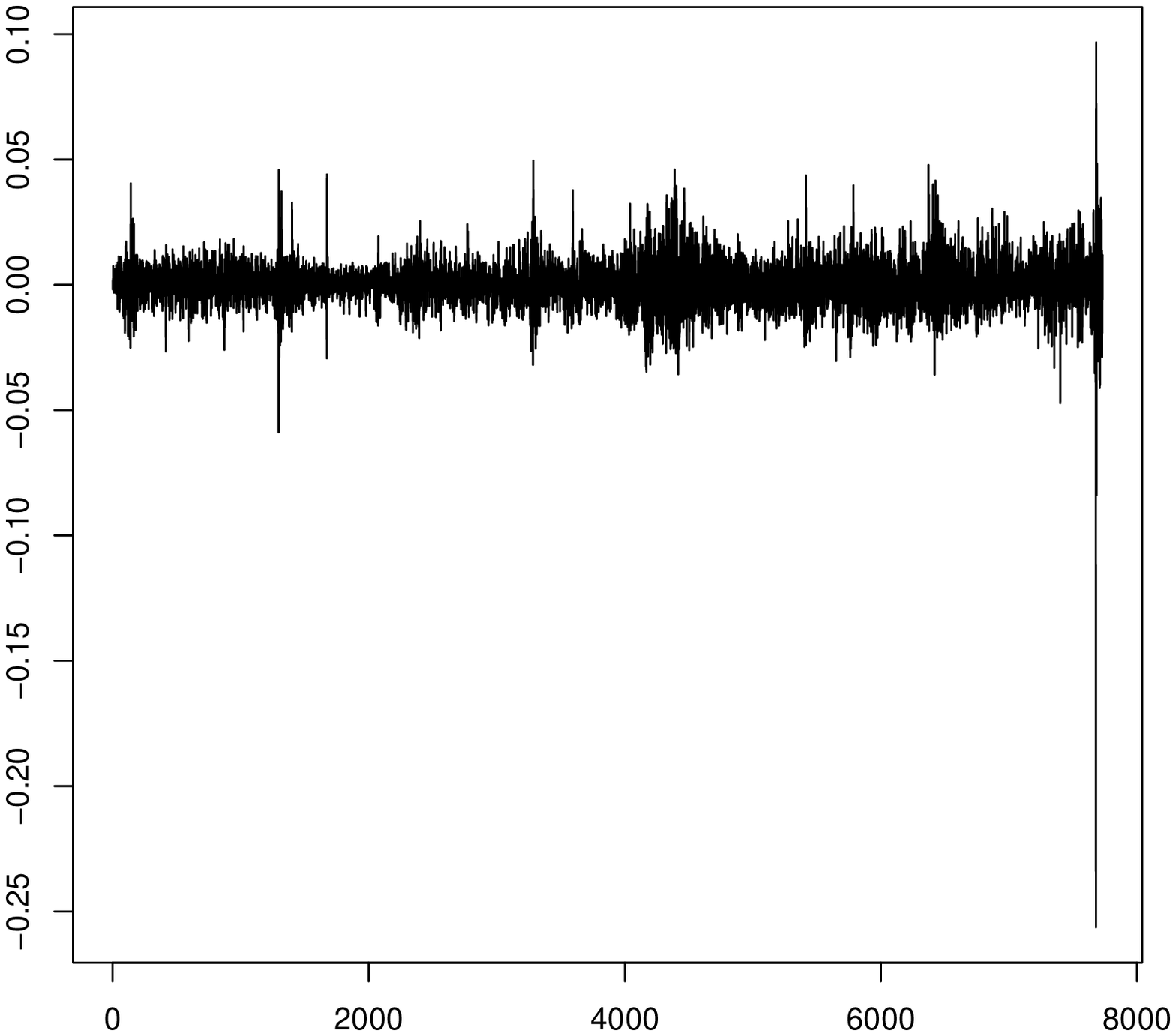}
\includegraphics[width= 5.8cm,height= 3.9cm]{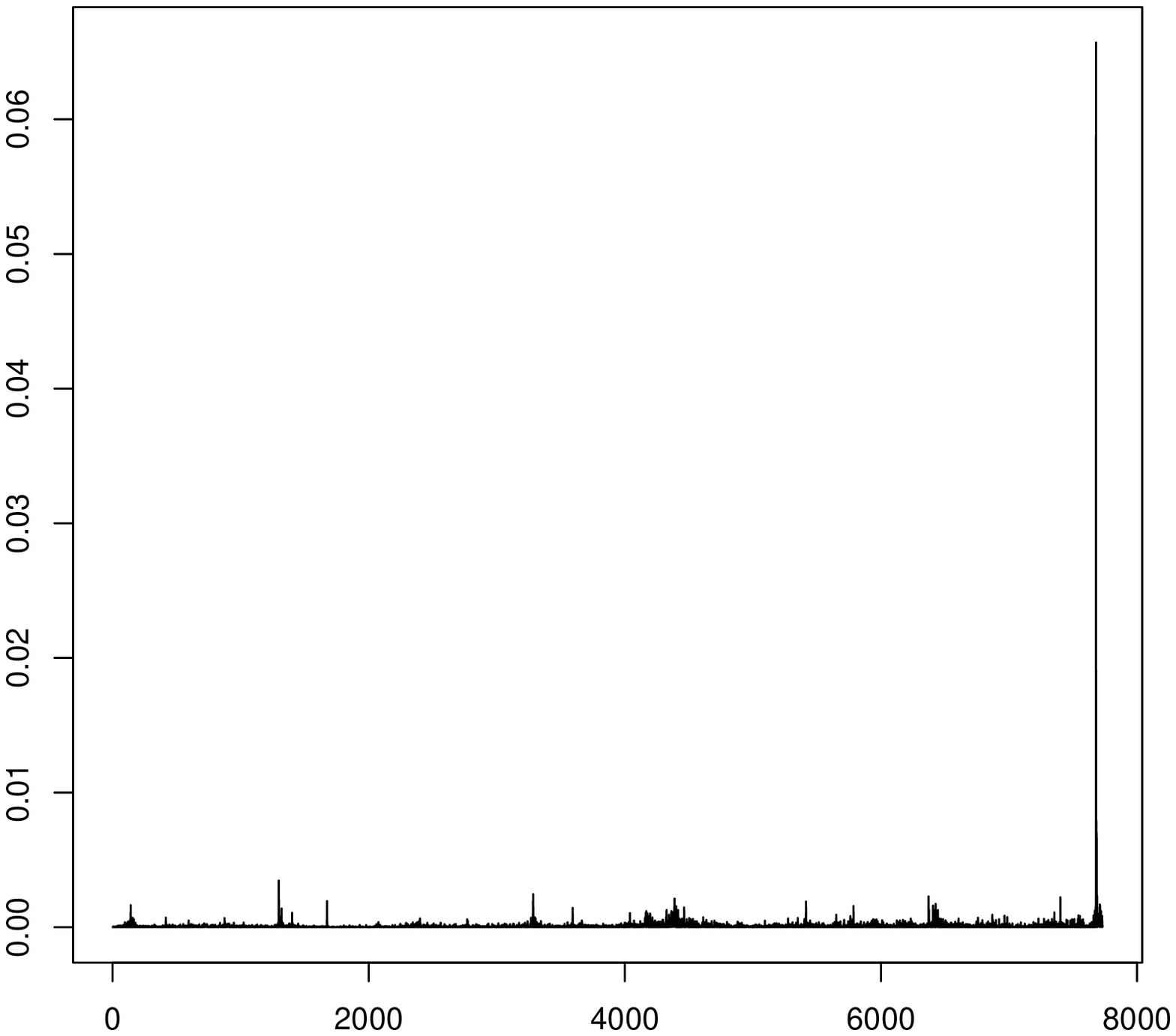}
\caption {Left: log-returns ($R_i$) of the S\&P500 (top) and Dow Jones indexes (bottom); Right: volatility ($R_i^2$) of the S\&P500 (top) and Dow Jones indexes (bottom).\label{figLR}}
\end{figure}

The presented results besides give us an insight of the effect of the random factors on the extremal behavior of the model, suggest estimation methods whenever data for $\mb{Y}$ and $\mb{T}$ are available, through the estimation of the dependence in $\mb{Y}$. Do not having real data sets in these conditions, we have chosen to illustrate in this section the application of the multivariate pRARMAX model considered in Example \ref{emrarmax} to a bivariate stock market index data. The univariate pRARMAX was introduced in Ferreira and Canto e Castro (\cite{lccmf2}, 2010) and an adjustment algorithm to real data was implemented and applied to the S\&P500 index, for the period April of 1957 to December 1987. Here we also consider the Dow Jones index for the same time period and our analysis is based on the volatility of the stock indexes measured through the square of the log-returns defined by $R_{i}=\log(P_{i+1})/\log(P_{i})$, $i=1,\hdots,7731$ ($n=7731$), where $P_{i}$ is the $i^{th}$ closing day of the index.  The values of $R_{i}^{\tt SP}$ and $\left(R_{i}^{\tt SP}\right)^2$, as well as, $R_{i}^{\tt DJ}$ and $\left(R_{i}^{\tt DJ}\right)^2$, corresponding, respectively, to the log-returns and volatility of the S\&P500 index and of the Dow Jones index, are plotted in Figure \ref{figLR}. Observe that they have similar plots. The large peak corresponds to ``Black Monday" stock market crash on the 19th October 1987. In Ferreira and Canto e Castro (\cite{lccmf2}, 2010) was considered a robust regression so that the marginals can be modeled through a Pareto distribution, and thus conducted the analysis on the transformed data $X^{SP}=a\,\left(R^{SP}\right)^2+b$, with $a=13618.3$ and $b=1.1$. Following the same procedure, we consider $X^{DJ}=a^*\,\left(R^{DJ}\right)^2+b^*$, with $a^*$ and $b^*$ estimated in $15154.2$ and $1.1$, respectively.

An  adjustment algorithm for a pRARMAX process, $X_n=Y_nZ_n$, $n\geq 1$, obtained by considering $Y_n=Y_{n-1}^c\vee W_n$ with $0<c<1$ and  $F_Y(x)=1-x^{-\beta}$, where $\{W_n\}_{n\geq 1}$ is an i.i.d. sequence independent of $Y_0$ and having distribution $F_W(x)=(1-x^{-\beta})/(1-x^{-\beta/c})$, and $\{Z_n\}_{n\geq 1}$ is a sequence of standard uniform random variables, was proposed in Ferreira and Canto e Castro (\cite{lccmf2} 2010, Section 3). This was used in this latter reference to model the $X^{SP}$ data producing the estimates $\widehat{c}=0.8$ (although the values $\widehat{c}=0.85,\,0.75$ were also considered) and $\widehat{\beta} =2$.
Now we are going to apply this adjustment algorithm to the $X^{DJ}$ data. In the following we summarize the running steps:
\begin{itemize}
\item[1.] Test if $X^{DJ}$  is in the Fr\'echet domain of attraction (see, e.g., Dietrich \emph{et al.} \cite{diet+} 2002) and estimate the tail index corresponding to $\beta^{-1}$ using, e.g., the Hill estimator (Hill, \cite{hill} 1975) and moments estimator (Dekkers \emph{et al.}, \cite{dek+} 1989).
\item[2.] Estimate the parameter c, through the estimation of the tail index of $S=\{S_i^{(n)}\}$, $i=1,\hdots, n-1$,  where $$S_i^{(n)}=\left(\frac{n+1}{n+1-K_n(X_i^{DJ})}\right)\wedge \left(\frac{n+1}{n+1-K_n(X_{i+1}^{DJ})}\right),$$ with $K_n(x)=n^{-1}\sum_{i=1}^n\mathds{1}_{\{X_i\leq x\}}$.
\item[3.] Apply the criterion: ``if $X_i^{DJ}>\left(X_{i-1}^{DJ}\right)^{\widehat{c}}$, with the estimate $\widehat{c}$ obtained in step 2, then $X_i^{DJ}$ comes from the innovation component $W_i$",
to separate the innovations $W$ and test if this
sample is also in the Fr\'echet domain of attraction.
\item[4.] Capture the observations corresponding to $Z$, through the criterion: ``if $X_i^{DJ}<\left(X_{i-1}^{DJ}\right)^{\widehat{c}}$ and $X_i^{DJ}\in \mathcal{B}_{\sscr\upsilon}=\left\{t:\frac{\pi_0 t^{-\widehat{\beta}/\widehat{c}}f_{W}(t)}{\pi_0 t^{ -\widehat{\beta}/\widehat{c}}f_{W}(t)+(1-\pi_0)t^{ -\widehat{\beta}/\widehat{c})-1}F_W (t)}\leq \upsilon\right\}$, with $f_W$ the density function of $W$, $\pi_1={P(Z_i\left(X_{i-1}^{DJ}\right)^{\widehat{c}}>W_i)}
    /{P(\left(X_{i-1}^{DJ}\right)^{\widehat{c}}>W_i)}$, $\pi_0=1-\pi_1$, and where $\widehat{\beta}$ and $\widehat{c}$ are the estimates obtained in steps 1.~and 2., respectively, then, $Z_i=X_i^{DJ}/\left(X_{i-1}^{DJ}\right)^{\widehat{c}}$ ''; $\upsilon$ cannot be too large (larger errors) nor too small (not enough observations to carry out the test of the next step).
\item[5.] Test whether the sample of random variables $Z$ captured in the step 4 has distribution $Beta(\widehat{\beta}/\widehat{c}+1,1)$ (e.g., Kolmogorov-Smirnov test).
\end{itemize}
For details, see Ferreira and Canto e Castro (\cite{lccmf2} 2010) and references therein.\\

\begin{figure}[!thb]
\centering
\includegraphics[width= 3.8cm,height= 3.9cm]{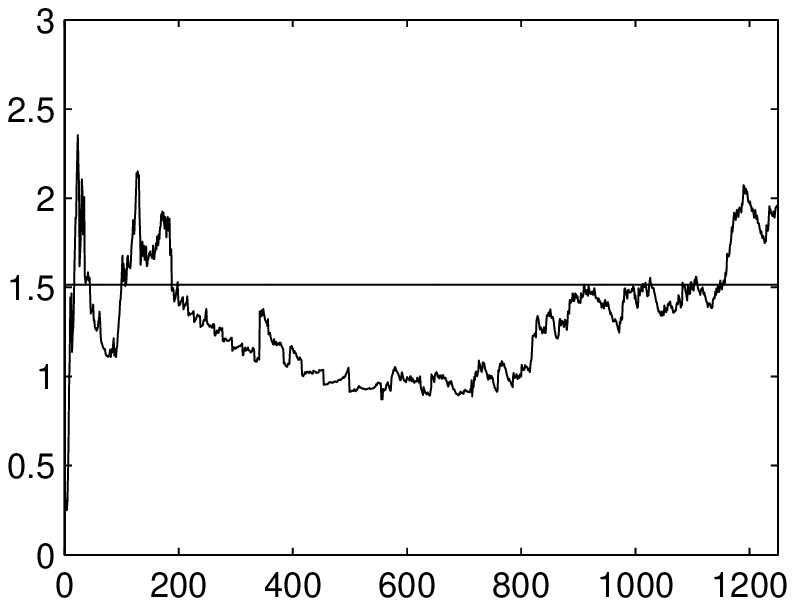}
\includegraphics[width= 3.8cm,height= 3.9cm]{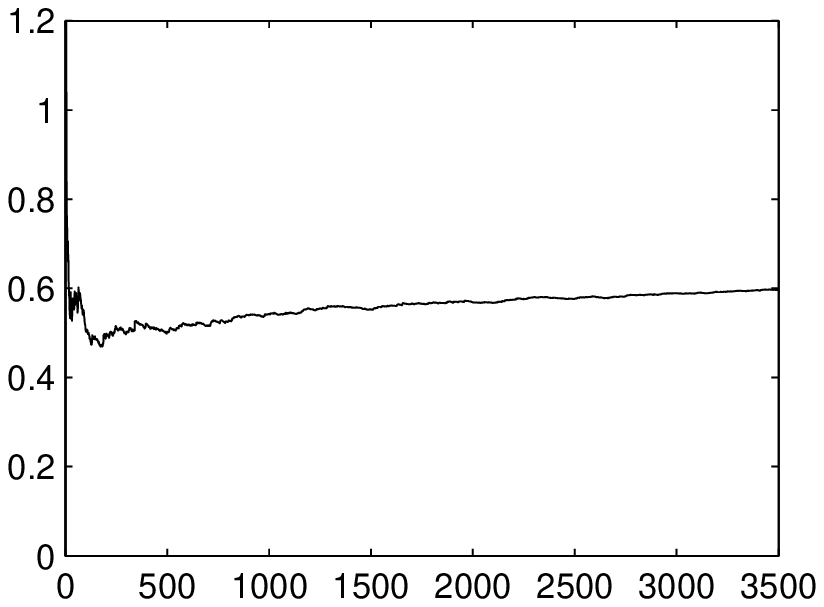}
\includegraphics[width= 3.8cm,height= 3.9cm]{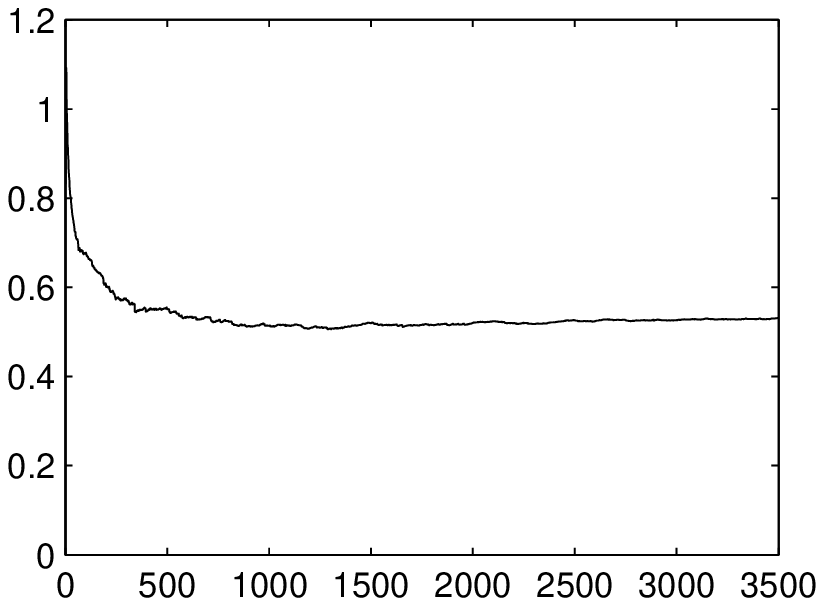}
\caption {Left: sample path of the heavy tailed test statistic, against k upper order statistics, applied to the data $X^{DJ}$ (horizontal line: critical value above which reject the null hypothesis of a heavy tailed domain of attraction); Center and right, respectively:
sample paths of Hill and moments estimators, against k upper order statistics, of the data $X^{DJ}$ concerning the estimation of the tail index $\beta^{-1}$.\label{figevtX}}
\end{figure}

In Figure \ref{figevtX} it is plotted the sample paths (against the $k$ upper order statistics) of the heavy tailed test statistic (left), the Hill estimator (center) and the moments estimator (right), referred in step 1. The heavy tailed domain of attraction is not rejected for $50\lesssim k\lesssim 120$ and $200\lesssim k\lesssim 1180$, which are plausible values within this text (see again Dietrich \emph{et al.} \cite{diet+} 2002). In what concerns
the sample paths of Hill and moments, the value where the paths yield approximately a flat line is at about $0.57$ and thus we consider $\widehat{\beta}=1/0.55\approx 1.82$.

\begin{figure}[!thb]
\centering
\includegraphics[width= 3.8cm,height= 3.9cm]{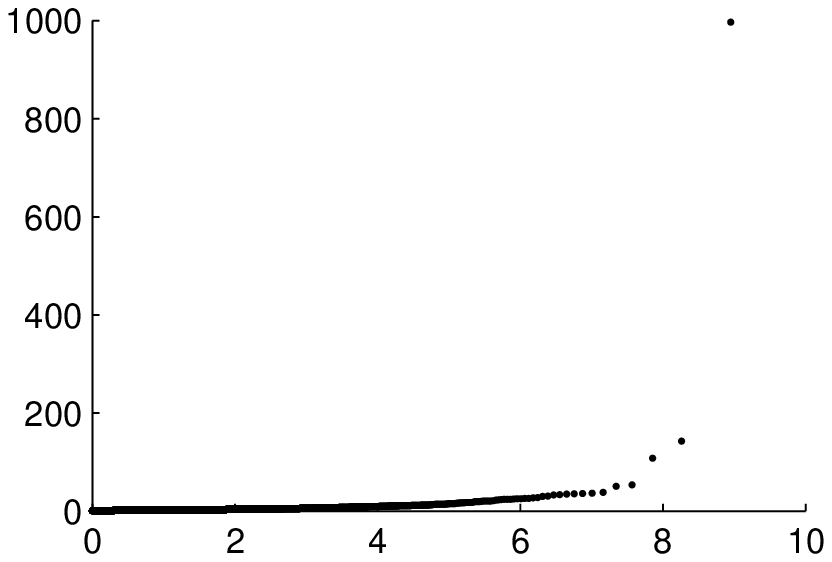}
\includegraphics[width= 3.8cm,height= 3.9cm]{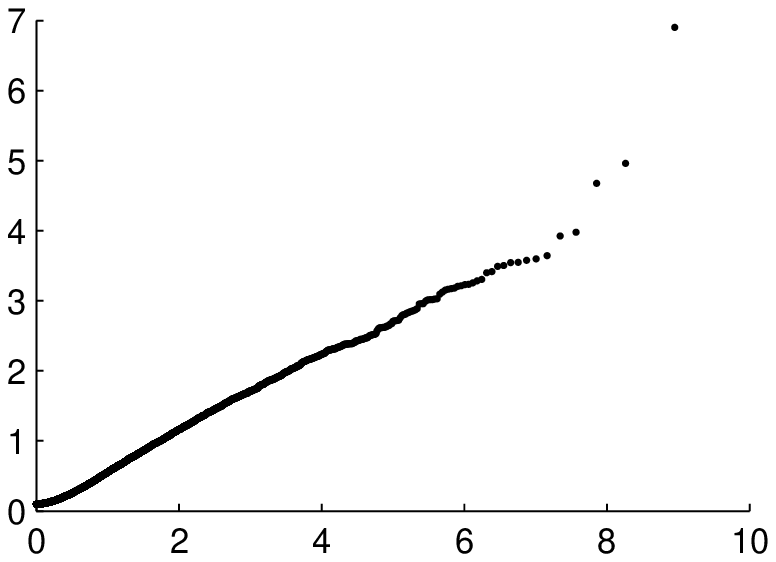}
\includegraphics[width= 3.8cm,height= 3.9cm]{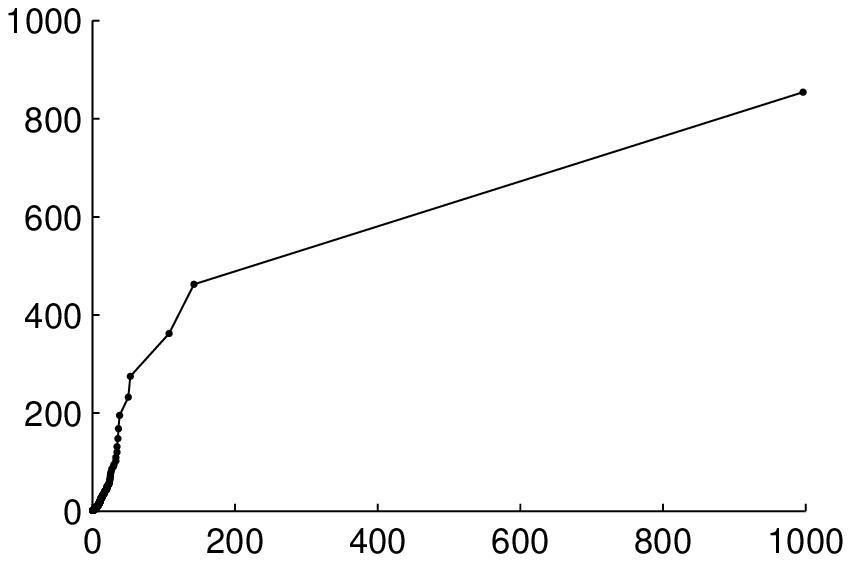}
\includegraphics[width= 3.8cm,height= 3.9cm]{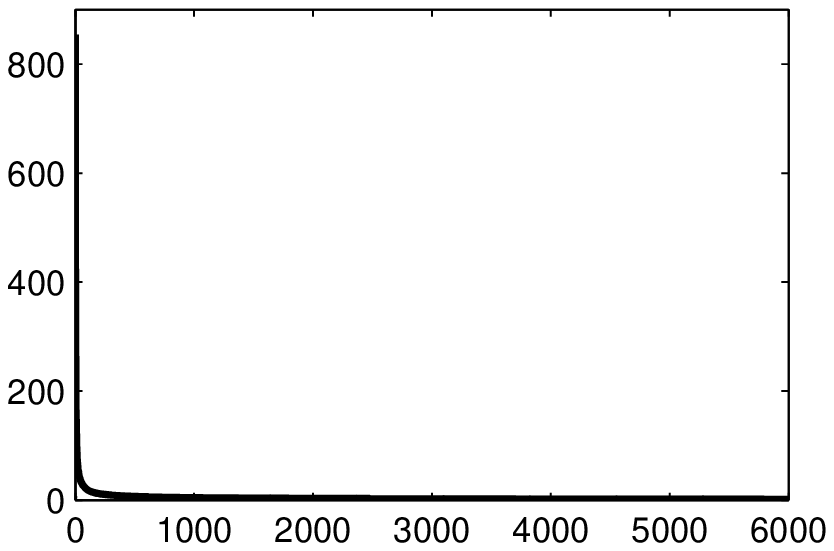}
\caption {Left to right, respectively: exponential and Pareto qq-plots of $X^{\tt DJ}$, empirical mean excess function against $X^{\tt DJ}$ and against k upper order statistics.\label{figqqs}}
\end{figure}

The Exponential and Pareto qq-plots and the empirical mean excess function plots (Figure \ref{figqqs}) indicate a Pareto-type model (see, e.g., Beirlant \emph{et al.} \cite{beirl+} 2004).

\begin{figure}[!thb]
\centering
\includegraphics[width= 3.8cm,height= 3.9cm]{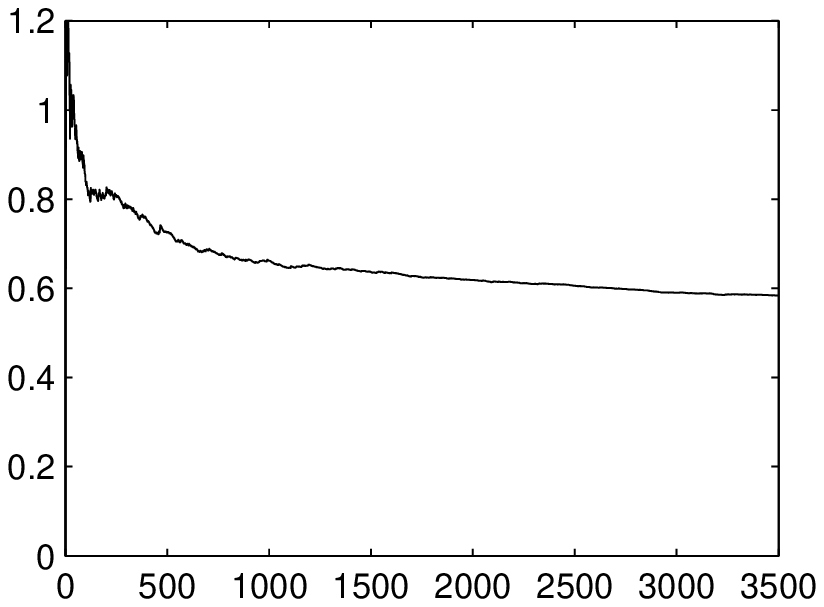}
\includegraphics[width= 3.8cm,height= 3.9cm]{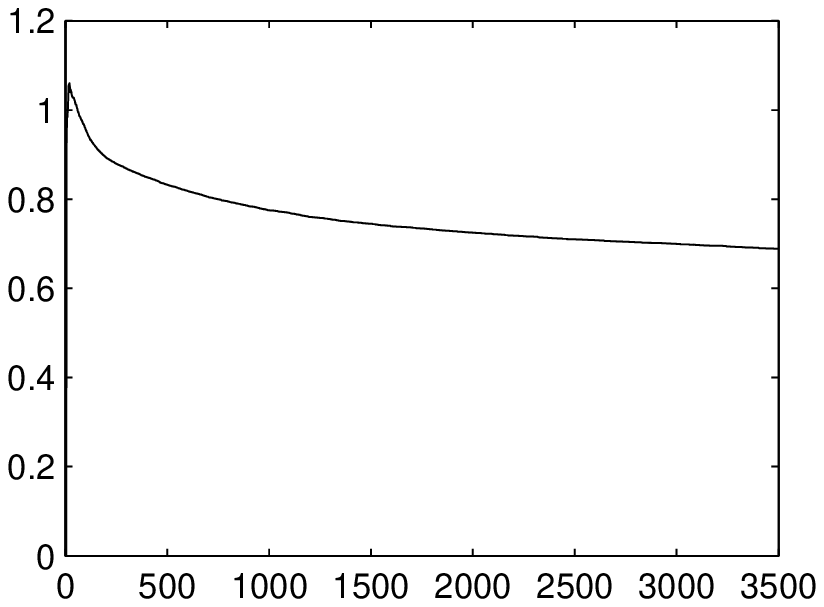}
\caption {Left and right, respectively:
sample paths of Hill and moments estimators, against k upper order statistics, of $S$ concerning the estimation of $c$.\label{figevtT}}
\end{figure}

In Figure \ref{figevtT} we have the sample paths of the Hill (left) and the moments (right) estimators concerning the sample $S$ (step 2). The plot is somewhat ``stable" at about $0.8$  and thus we take $\widehat{c}=0.8$.

\begin{figure}[!thb]
\centering
\includegraphics[width= 3.8cm,height= 3.9cm]{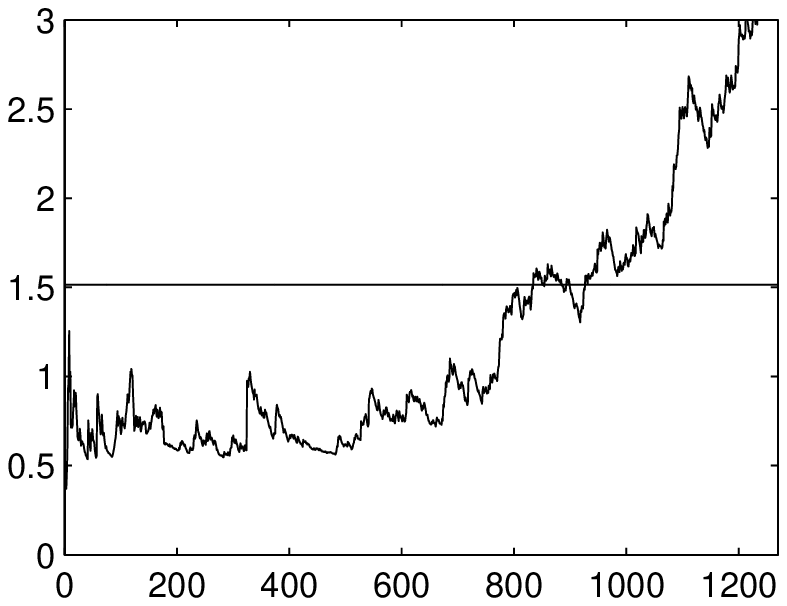}
\includegraphics[width= 3.8cm,height= 3.9cm]{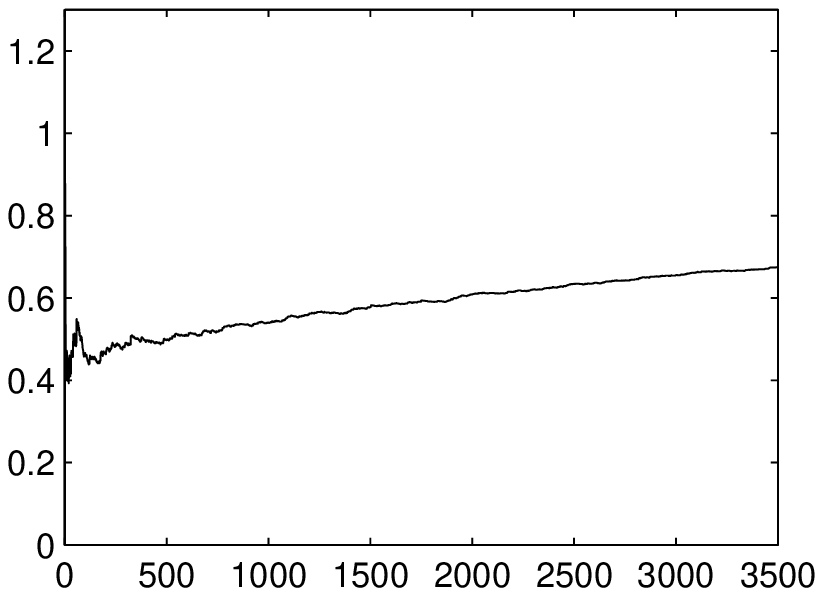}
\includegraphics[width= 3.8cm,height= 3.9cm]{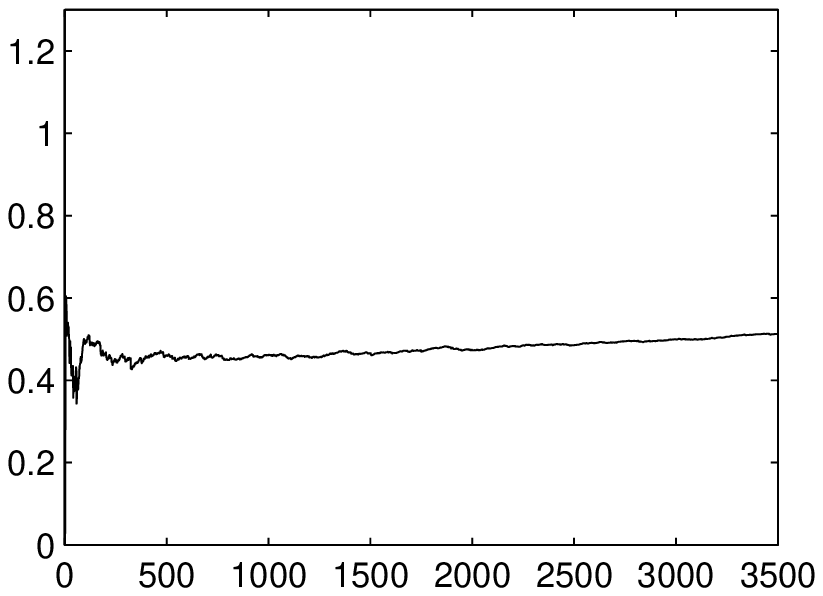}
\caption {Left: sample path of the heavy tailed test statistic, against k upper order statistics, applied to the separated innovations $W$, according to step 3 (horizontal line: critical value above which reject the null hypothesis of a heavy tailed domain of attraction); Center and right, respectively:
sample paths of Hill and moments estimators, against k upper order statistics, of W concerning the estimation of the tail index $\beta^{-1}$.\label{figevtZ}}
\end{figure}

Applying the criterium of step 3, we separate the innovations component $W$. The sample path of the heavy tailed test statistic (Figure \ref{figevtZ}, left) do not reject this assumption for $k\lesssim 780$, which is again a plausible result to conclude a Fr\'echet domain of attraction. From the Hill and moments plots (Figure \ref{figevtZ}, center and right, respectively) we do not discard the same estimate of $\widehat{\beta}\approx 1.82$.

According to step 4, we capture the observations corresponding to random coefficients, $Z$, considering the significance regions $\mathcal{B}_\upsilon$ for $\upsilon= 0.05,0.1,\hdots,0.45$ and $\widehat{c}=0.8$.

\begin{figure}[!thb]
\centering
\includegraphics[width= 3.8cm,height= 3.9cm]{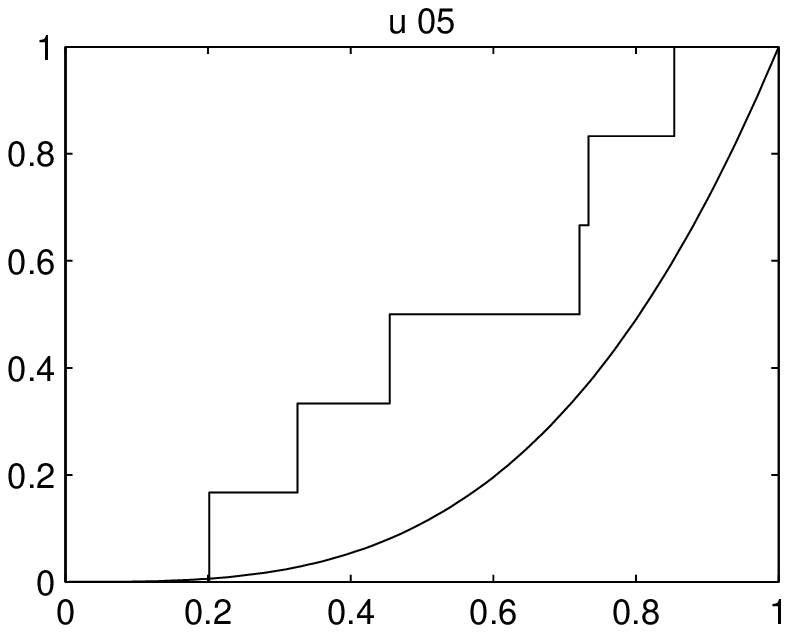}
\includegraphics[width= 3.8cm,height= 3.9cm]{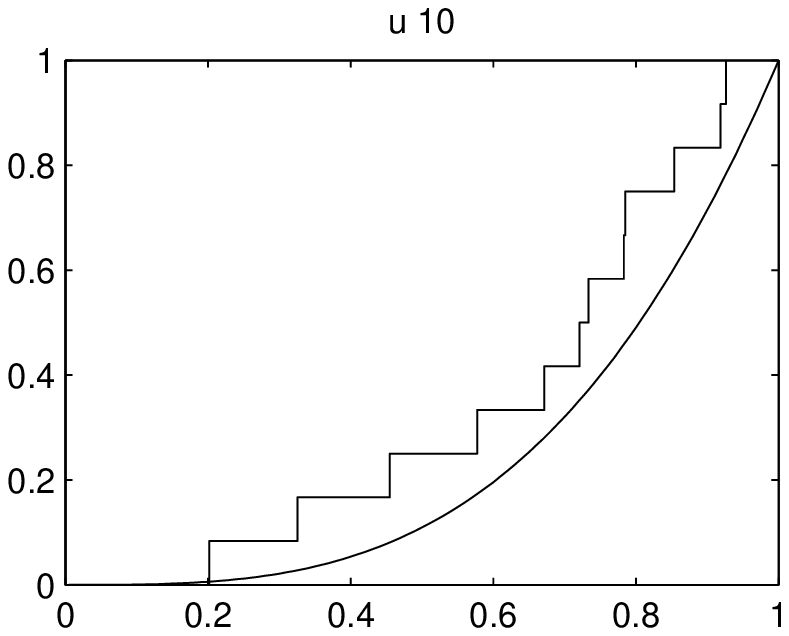}
\includegraphics[width= 3.8cm,height= 3.9cm]{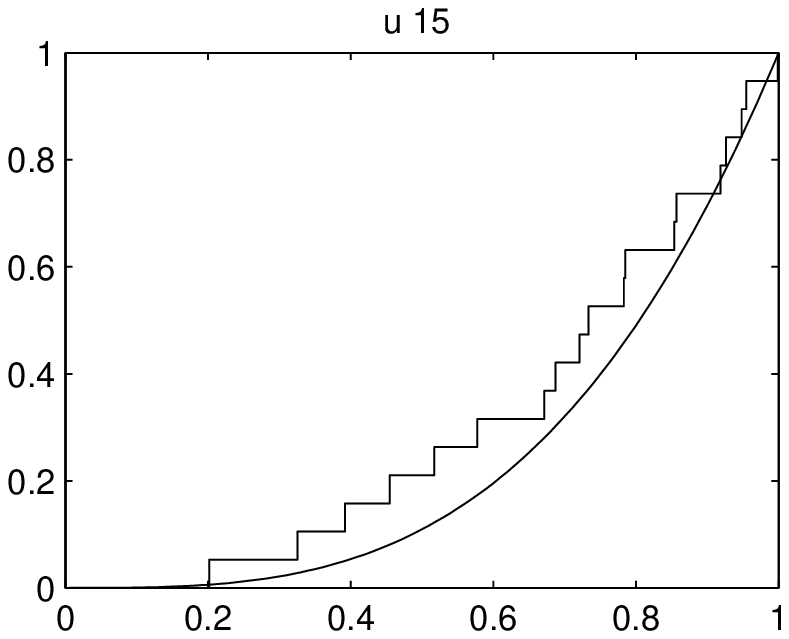}\\
\includegraphics[width= 3.8cm,height= 3.9cm]{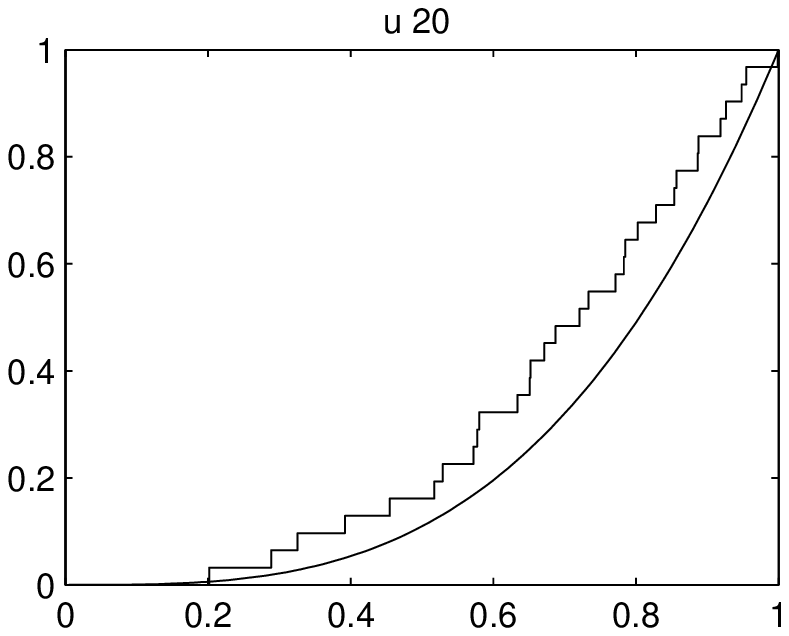}
\includegraphics[width= 3.8cm,height= 3.9cm]{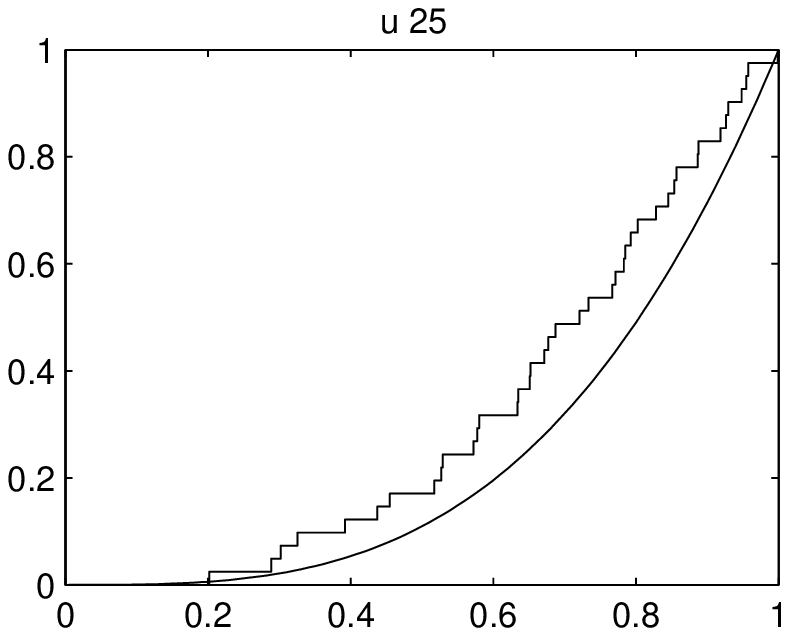}
\includegraphics[width= 3.8cm,height= 3.9cm]{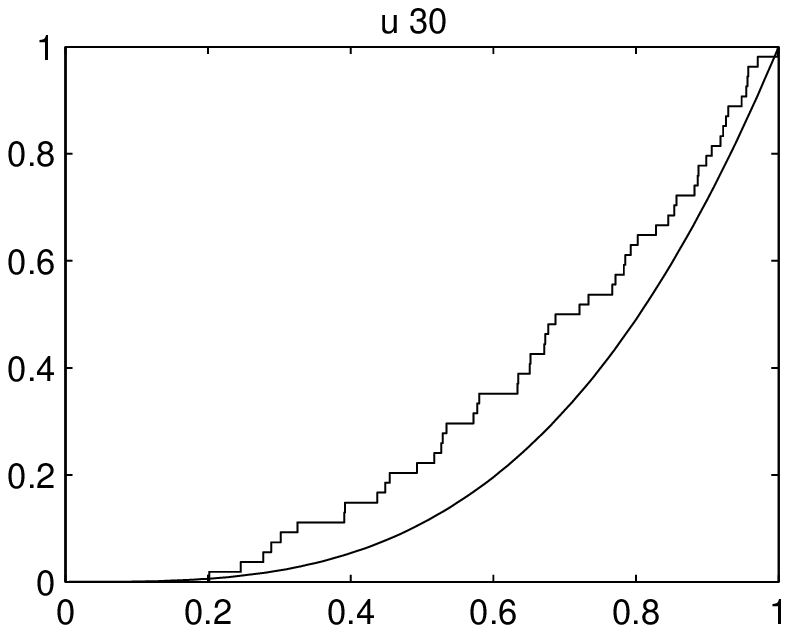}\\
\includegraphics[width= 3.8cm,height= 3.9cm]{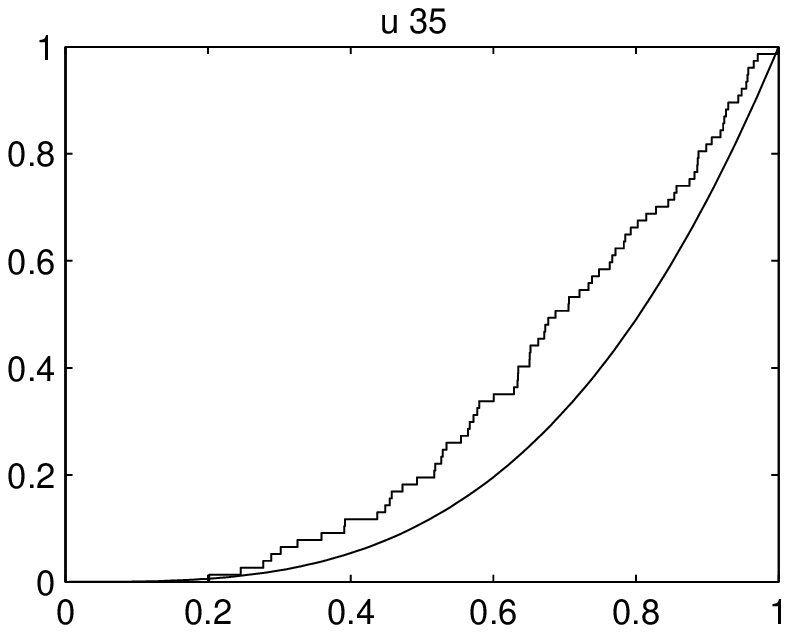}
\includegraphics[width= 3.8cm,height= 3.9cm]{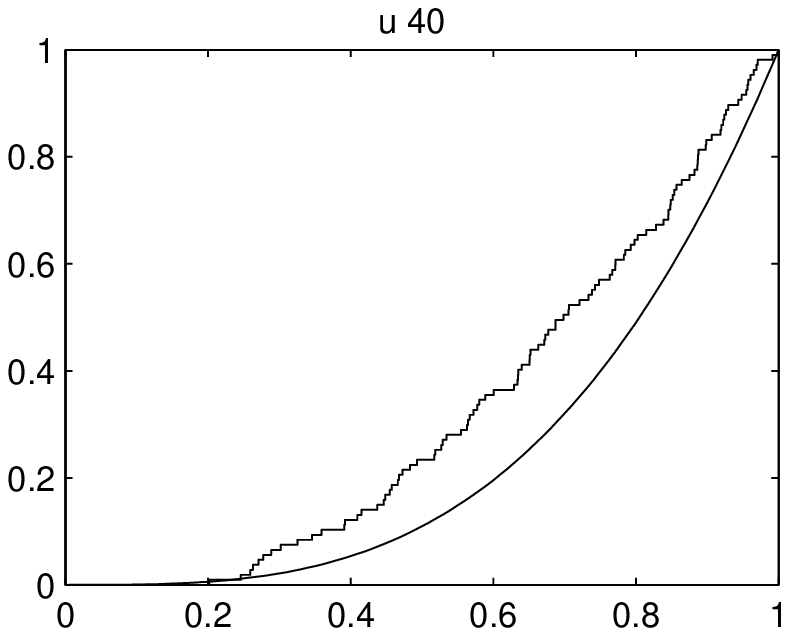}
\includegraphics[width= 3.8cm,height= 3.9cm]{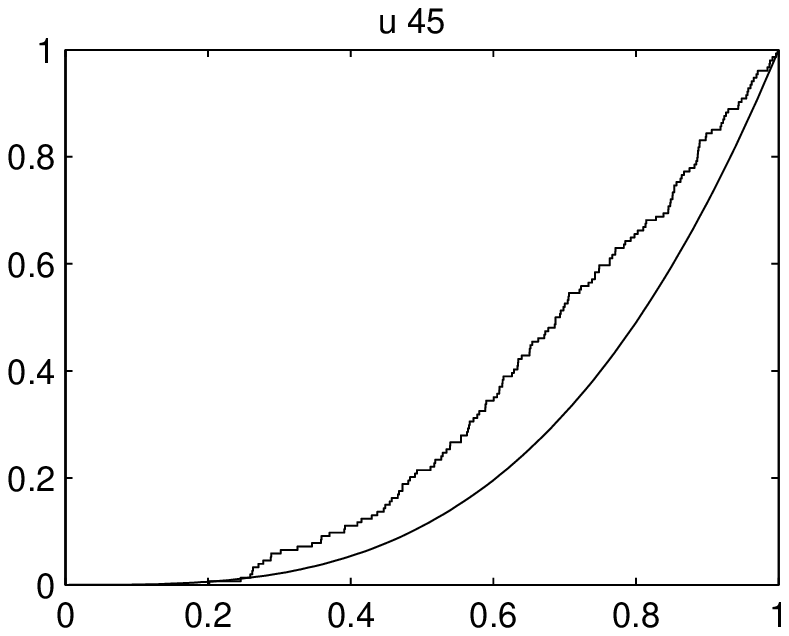}
\caption {Empirical and theoretical distribution functions of the random coefficients $Z$ captured on step 4, for significance regions with,
$\upsilon=0.05$ (``u 05"),$ \hdots$, $\upsilon=0.45$ (``u 45").\label{figks}}
\end{figure}

Implementing step 5, we applied the Kolmogorov-Smirnov test for the distribution ${\tt Beta}(1.82/0.8+1,1)$ to the sample of random variables $Z$ captured in the step 4. The rejection is obtained for $\upsilon\geq 0.3$. By Figure \ref{figks}, the value $\upsilon=0.15$ seems to be an appropriate choice (with $21$ observations captured) and is in agreement with the simulation study conducted in Ferreira and Canto e Castro (\cite{lccmf2} 2010, Section 4.1). Therefore, we consider that the model pRARMAX with $c=0.8$ and $\beta=1.82$, can be used for the modeling of the transformed data $X^{\tt DJ}$.

\begin{figure}[!thb]
\centering
\includegraphics[width= 4.8cm,height= 4.9cm]{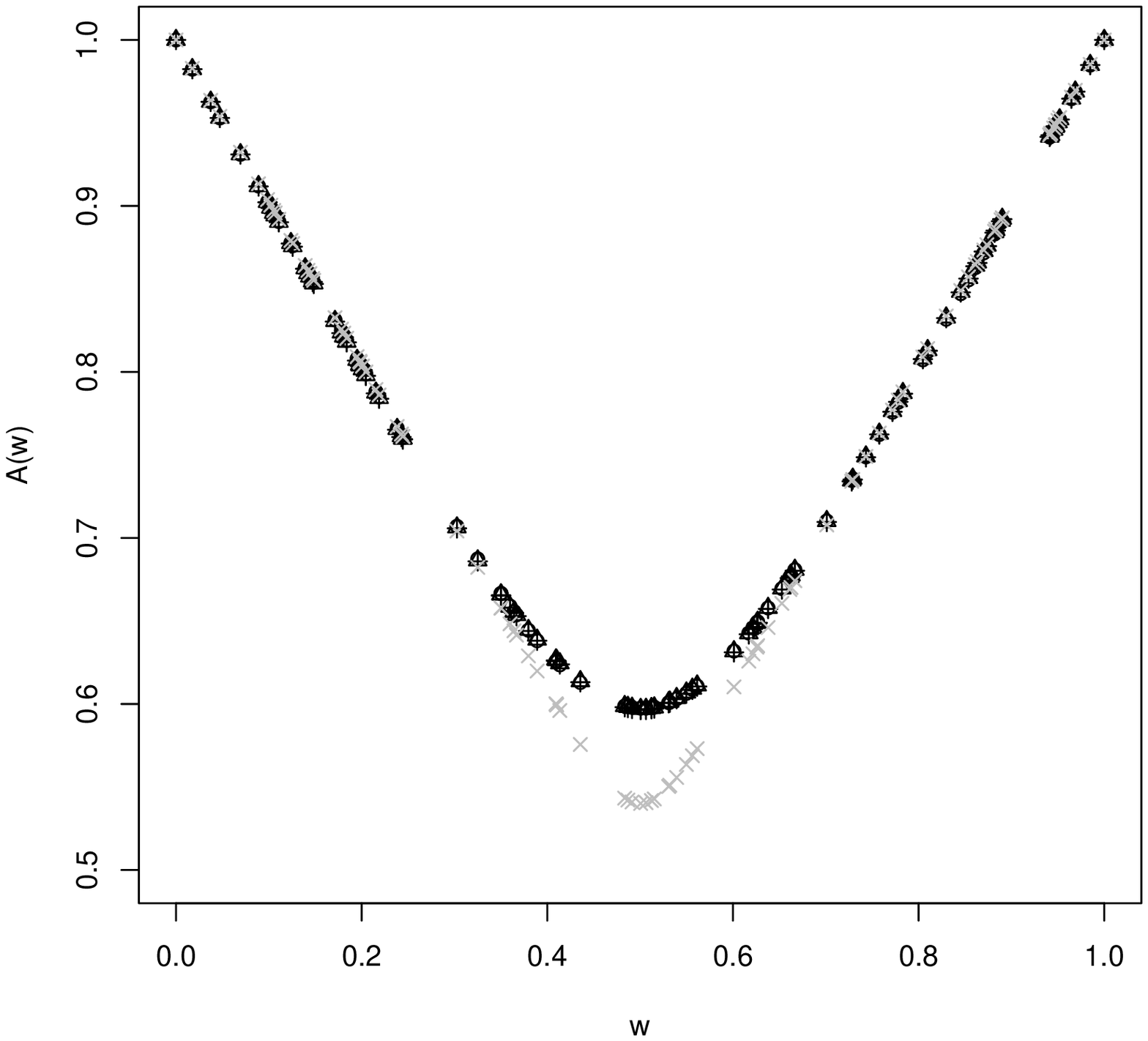}
\includegraphics[width= 4.8cm,height= 4.9cm]{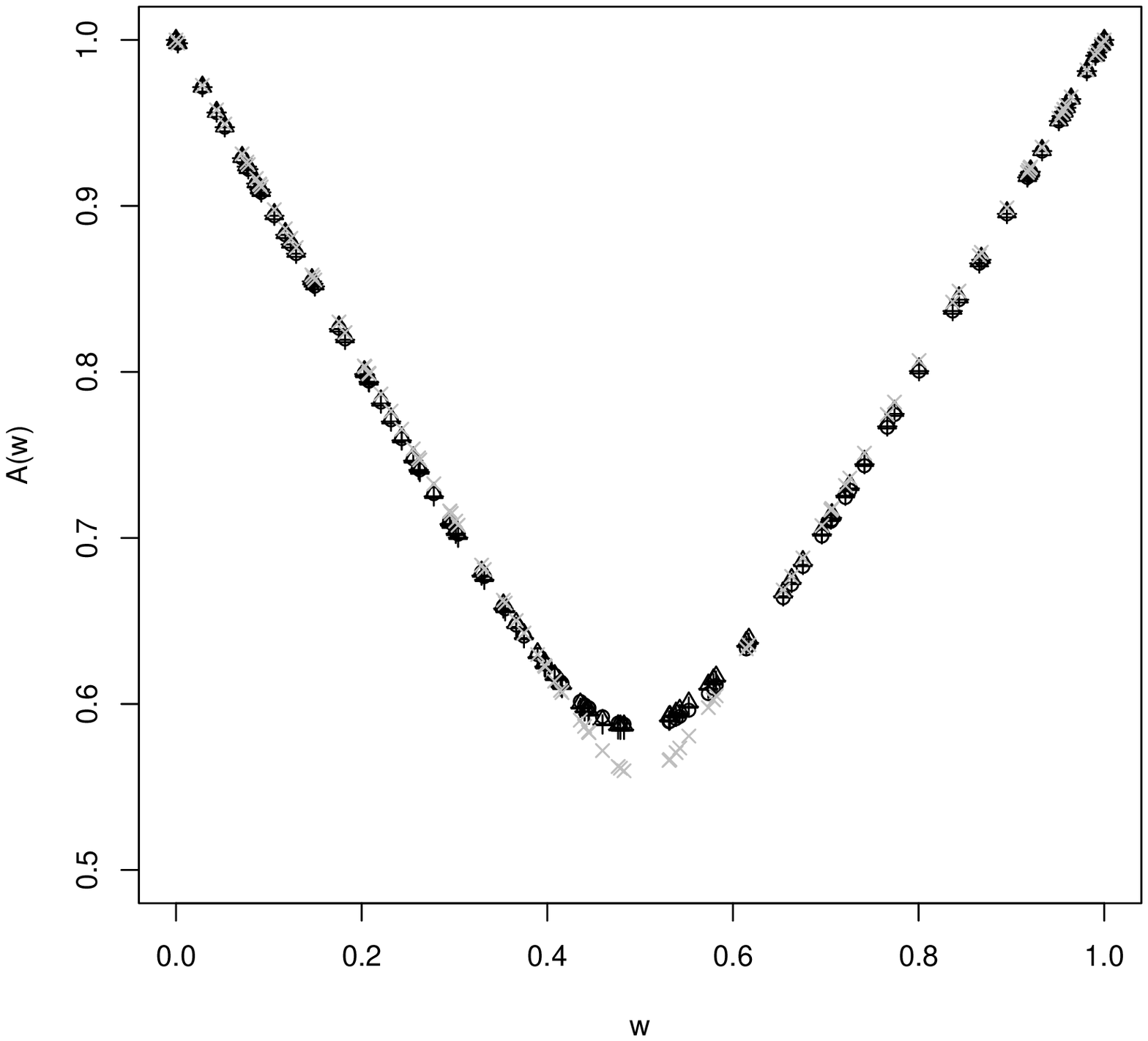}
\includegraphics[width= 4.8cm,height= 4.9cm]{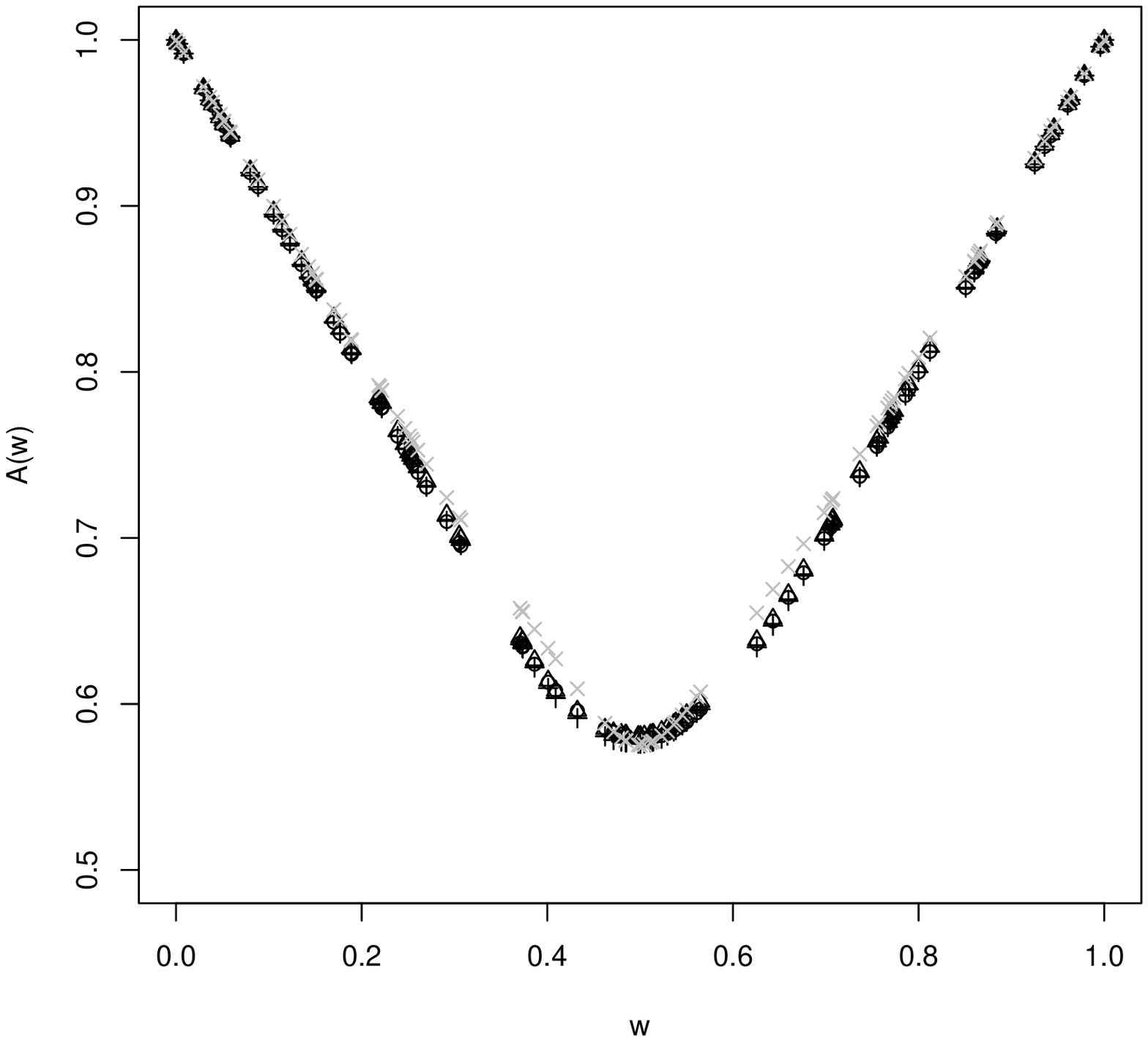}\\
\includegraphics[width= 4.8cm,height= 4.9cm]{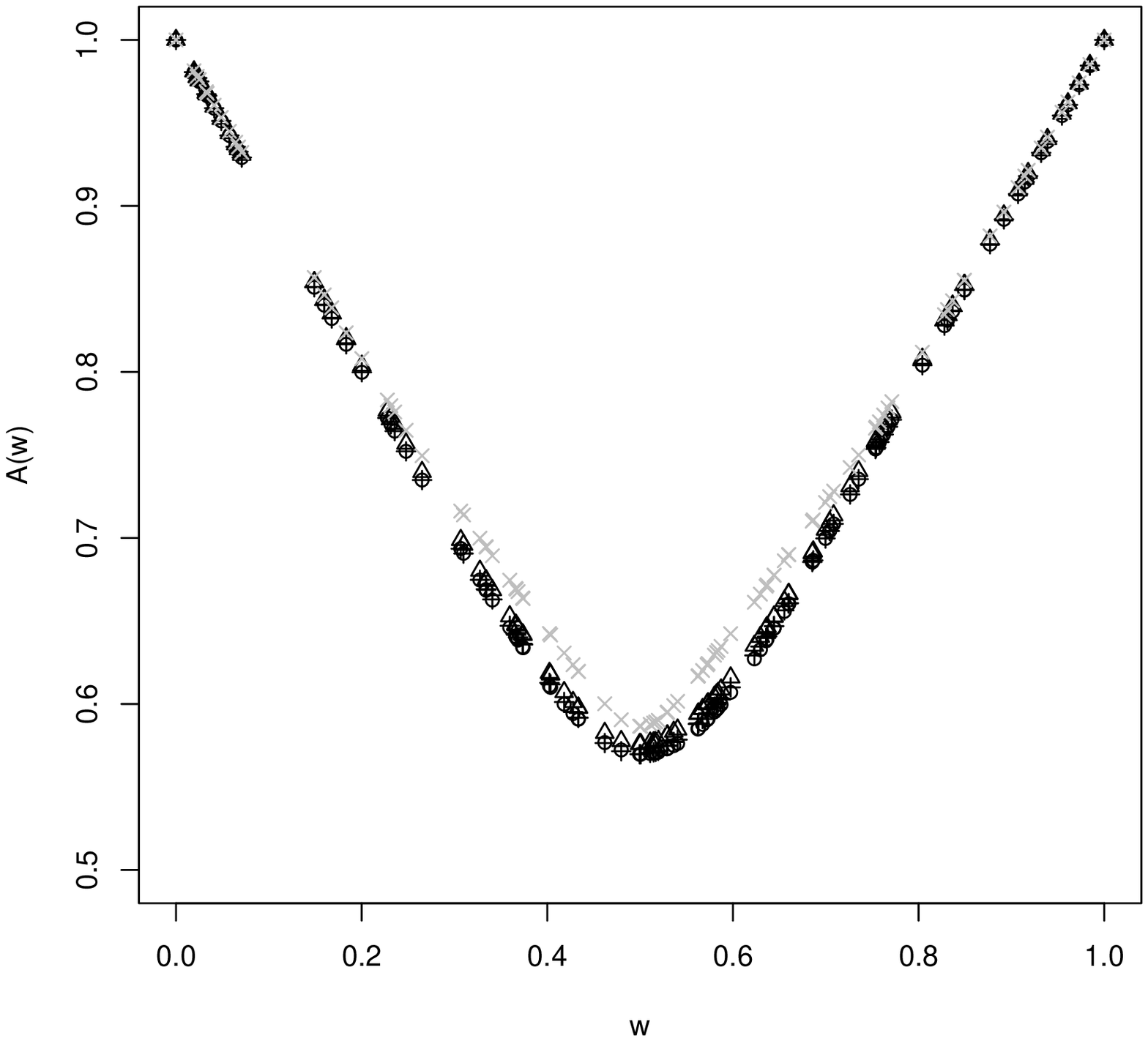}
\includegraphics[width= 4.8cm,height= 4.9cm]{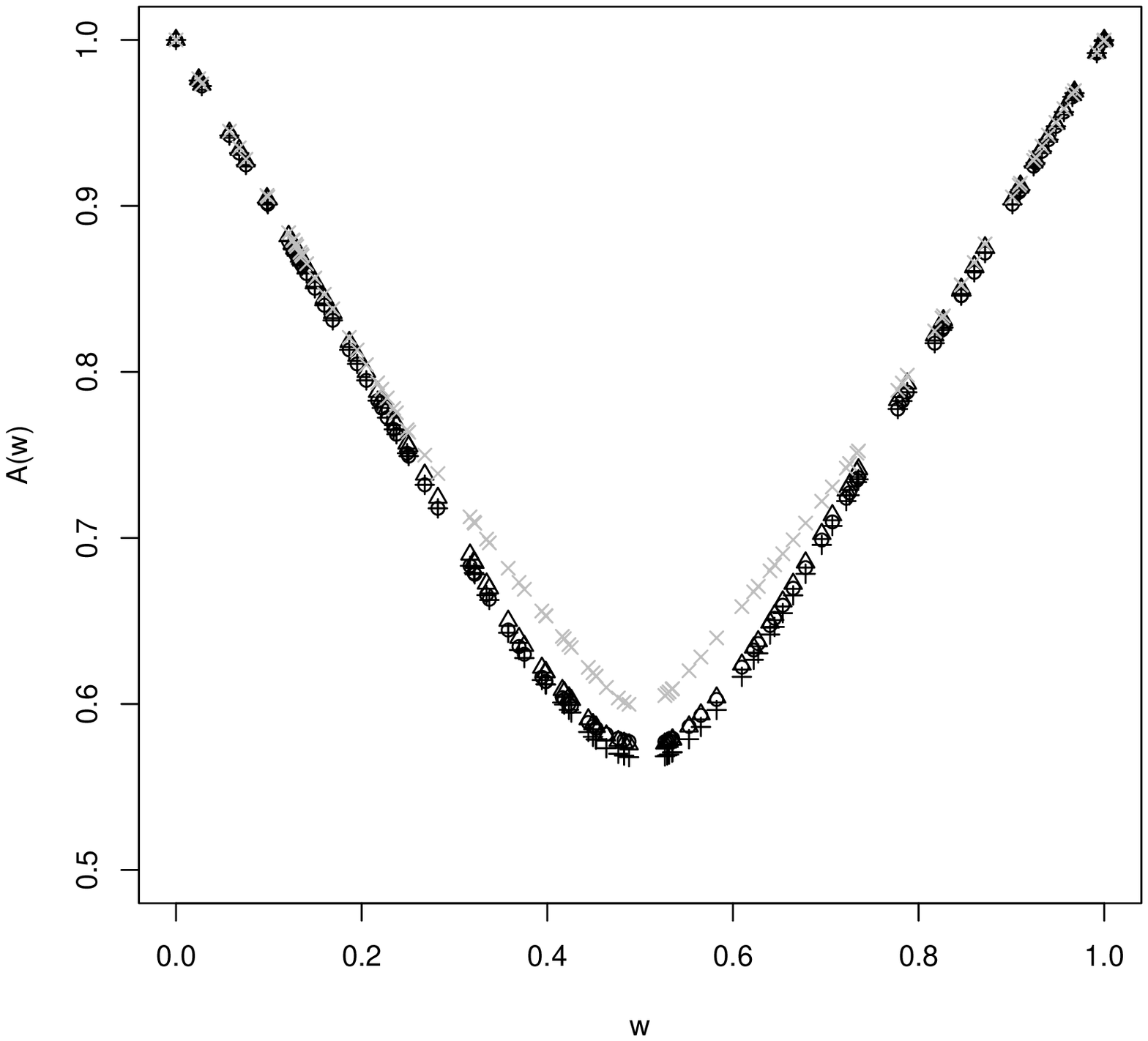}
\includegraphics[width= 4.8cm,height= 4.9cm]{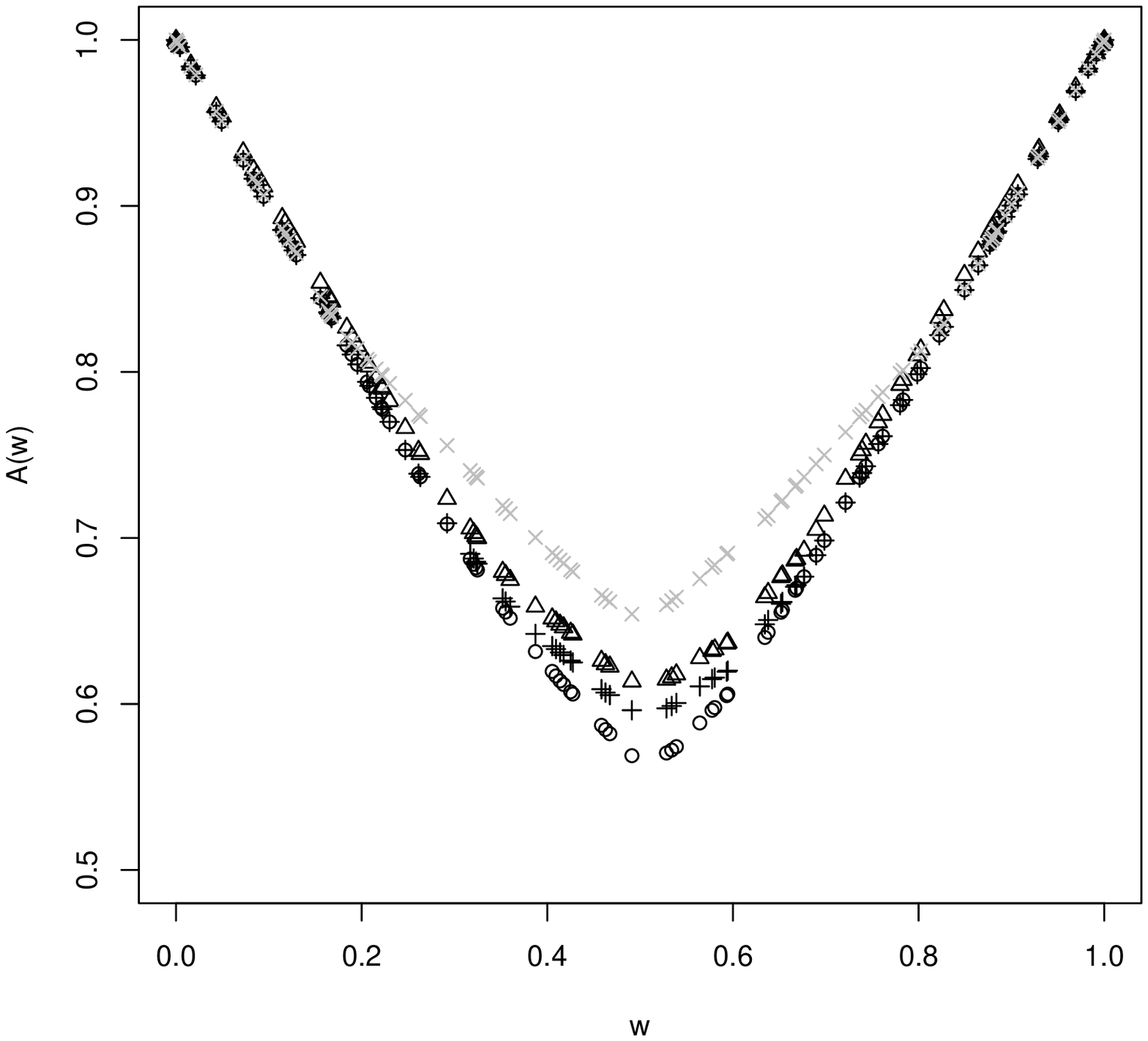}\\
\caption {Estimation of the dependence function $A$ in (\ref{depA}) for the bivariate data $\left({X_i}^{\tt SP},{X_i}^{\tt DJ}\right)$, $i=1,\hdots,7731$, using the non parametric estimators of Cap\'era\`a \emph{et al.} \cite{cap+97} ($\circ$), Pickands \cite{pic81} ($+$), Hall and
Tajvidi \cite{hall+taj04} ($\vartriangle$) and using the estimator of Lescourret and Robert \cite{les+rob06} ($\times$; grey), considering (top-to-bottom and left-to-right, respectively) the weekly, monthly, bimonthly, three monthly, quarterly and anual maxima.\label{figdep}}
\end{figure}

Now we are going to estimate the dependence of $(X_i^{\tt SP},X_i^{\tt DJ})$, $i=1,\hdots,7731$, based on the  dependence function $A$ (Pickands \cite{pic81} 1981) of a GEV distribution $G$, defined by
$$
G(x_1,x_2)=\exp\left\{-(x_1^{-1}+x_2^{-1})A(x_1(x_1+x_2)^{-1})\right\}.
$$
In order to have a dataset more coherent with an MEV model, we take block maxima within each marginal. We use the Pickands (\cite{pic81} 1981), the Cap\'er\`a \emph{et al.} (\cite{cap+97} 1997) and the Hall and Tajvidi (\cite{hall+taj04} 2004) nonparametric estimators for the function $A$, and thus free of model assumptions underlying the data. We also consider the Lescourret and Robert estimator (\cite{les+rob06}, 2006), given by
\begin{eqnarray}\label{depA}
\widehat{A}(w)=\min\left(1,\max\left(\frac{\sum_{i=1}^n \max\left((1-w)\widetilde{X_i}^{\tt SP},w\widetilde{X_i}^{\tt DJ}\right)}{\sum_{i=1}^n \left((1-w)\widetilde{X_i}^{\tt SP}+w\widetilde{X_i}^{\tt DJ}\right)},w,1-w\right)\right),
\end{eqnarray}
with
$$
\left(\widetilde{X_i}^{\tt SP},\widetilde{X_i}^{\tt DJ}\right)= \left(\left(\frac{{X_i}^{\tt SP}}{{X_{n-k:n}}^{\tt SP}}\right)^{\widehat{\beta}_{SP}},\left(\frac{{X_i}^{\tt DJ}}{{X_{n-k:n}}^{\tt DJ}}\right)^{\widehat{\beta}_{DJ}}\right)\,\mathds{1}_{\{{X_i}^{\tt SP}>{X_{n-k:n}}^{\tt SP},{X_i}^{\tt DJ}>{X_{n-k:n}}^{\tt DJ}\}},
$$
where
${X_{1:n}}^{\tt j}\leq \hdots \leq {X_{n:n}}^{\tt j}$, $j\in\{\tt SP,DJ\}$, are the respective order statistics. We have considered $\widehat{A}(w)=(\widehat{A}(w)+\widehat{A}(1-w))/2$ in order to constrain the estimators to be symmetric.  This estimator was developed for factor models with a common factor  $Y$, as described in Section \ref{spc}. The multivariate pRARMAX model  (Example \ref{emrarmax}) was derived within this context. Observe in Figure \ref{figdep} that the estimators are quite close, specially if we consider the bimonthly maxima, an indication that the multivariate pRARMAX may be a plausible model.

\end{document}